%% file: main.tex
\documentclass[11pt,a4paper]{article}

\usepackage[english]{babel}
\usepackage[utf8]{inputenc}
\DeclareUnicodeCharacter{2212}{-}
\usepackage{amsmath,amssymb,amsbsy,amsfonts}
\usepackage{bm}
\usepackage{stmaryrd}

\usepackage{ragged2e}
\justifying

\usepackage[paperwidth=192mm,
		    paperheight=262mm,
		    vmargin={19mm,19mm},
		    hmargin={13.7mm,13.7mm},
		    headsep=12pt,
		    footskip=12pt]{geometry}

\usepackage{hyperref}

\usepackage{siunitx}
\usepackage{graphicx}
\usepackage{verbatim}
\graphicspath{{./figures/}}
\usepackage{float}
\usepackage{caption}
\usepackage{subcaption}

\usepackage{tabularx}
\usepackage{multirow}
\usepackage{booktabs}

\usepackage{csquotes}

\usepackage[
  backend=biber,
  style=ieee,
  sorting=none,
  doi=false,
  url=false,
  isbn=false,
  eprint=false,
  giveninits=true,
  date=year,
  citestyle=ieee-comp
]{biblatex}

\AtEveryBibitem{
  \clearlist{language}
  \clearfield{month}
  \clearfield{day}
  \clearfield{series}
  \clearfield{number}
  \clearfield{urlyear}
  \clearfield{urlmonth}
  \clearfield{urlday}
}

\renewbibmacro{in:}{}

\DeclareFieldFormat
  [article,inproceedings,incollection, misc]
  {title}{#1}

\newcommand{\vel}{\bm{u}}

\newcommand{\jump}[1]{\llbracket{#1}\rrbracket}
\newcommand{\tensorjump}[1]{\langle\!\langle {#1} \rangle\!\rangle}
\newcommand{\avg}[1]{\{\!\!\{{#1}\}\!\!\}}

\begin{document}

\title{A high-order IMEX-DG solver for atmospheric dynamics: rotating tests on non-conforming meshes}

\author{Letizia Bottani$^{(1,2)\,,\,a}$, Tommaso Benacchio$^{(3)}$, Giuseppe Orlando$^{(4),\,b}$, \\ Luca Bonaventura$^{(1)}$, Allan Peter Engsig-Karup$^{(2)}$}

\maketitle

\begin{center}
{
    \small
   $^{(1)}$
   Dipartimento di Matematica, Politecnico di Milano \\
    Piazza Leonardo da Vinci 32, 20133 Milano, Italy \\
    %
    \vspace{3mm}
    $^{(2)}$
    Department of Applied Mathematics and Computer Science, Technical University of Denmark, Asmussens Allé, 303B, Kongens Lyngby, 2800, Denmark \\
    %
    \vspace{3mm}
    $^{(3)}$
    Weather Research, Danish Meteorological Institute \\
    Sankt Kjelds Plads 11, 2100 Copenhagen, Denmark \\
    
    %
    \vspace{3mm}
    $^{(4)}$
    CMAP, CNRS, \'{E}cole polytechnique, Institut Polytechnique de Paris \\
    Route de Saclay, 91120 Palaiseau, France \\
    \vspace{3mm}
    $^{(a)}$ Now at Chair of Computational Methods for Environmental and Civil Engineering, Technical University of Munich, Theresienstraße 90, 80333
Munich, Germany\\
 \vspace{3mm}
$^{(b)}$ Now at Institut de Mathématiques de Toulouse, UMR CNRS 5219, Université de Toulouse, INSA, F-31077 Toulouse, France.
\vspace{3mm}\\
Email adresses: {\tt letizia.bottani@tum.de, tbo@dmi.dk, giuseppe.orlando@math.univ-toulouse.fr, luca.bonaventura@polimi.it , apek@dtu.dk}
 
}
\end{center}

\vspace{3mm}


\pagebreak

\begin{abstract}

We present a high-order implicit-explicit discontinuous Galerkin
(IMEX-DG) solver for the compressible Euler equations that accounts
for rotational effects within an atmospheric
framework. Time integration follows a second-order additive
Runge--Kutta scheme, treating stiff acoustic modes implicitly and
advective terms explicitly. The solver is built on the
\texttt{deal.II} finite element library, combining matrix-free
operator evaluation, adaptive non-conforming mesh capabilities,
and distributed-memory parallelism.

The inclusion of rotation extends the framework to large-scale
atmospheric dynamics, where the Coriolis force plays a central role
in geostrophic balance. Its incorporation into a high-order
matrix-free DG discretisation is, however, non-trivial, as it
requires a careful formulation of the discrete rotational operator.
Two alternative treatments of the rotational and gravitational
source terms within the solution strategy, based on nonlinear
fixed-point iterations, are therefore introduced and compared in
terms of accuracy, robustness, and computational efficiency.
A discrete analysis of the rotational operator is also carried out
in order to derive a formulation suitable for efficient matrix-free
implementation.

The proposed formulation is verified through convergence studies
on rotating inertia-gravity wave benchmarks and further assessed
in fully three-dimensional simulations of stratified flow over
orography on both uniform and adaptive meshes. The numerical
results show that the rotating IMEX-DG framework has the expected
accuracy and stability properties while correctly capturing the
asymmetry and wave structures induced by rotation in large-scale
atmospheric flows.
\end{abstract}
\setlength{\leftskip}{0pt}
\setlength{\rightskip}{0pt}

\pagebreak

\input{Introduction}

\input{Math_model}

\input{Numerical_Framework}

\input{Numerical_results}

\input{Conclusion}


\printbibliography
\end{document}

%% file: Introduction.tex
\section{Introduction}
\label{sec:intro}

Atmospheric flows are characterized by the interaction of multiple spatial and temporal scales, ranging from large-scale balanced circulations to fast acoustic and gravity waves \cite{holton_introduction_2013,vallis_atmospheric_2017}. For this reason, fully compressible formulations are widely used in atmospheric modeling, as they provide a complete description of these dynamics \cite{steppeler_review_2003}. In particular, although atmospheric flows typically exhibit relatively low flow speeds, compressibility effects cannot in general be neglected because of strong vertical stratification and associated gravity-wave dynamics. However, in fully compressible atmospheric models, the coexistence of slow geostrophic modes and advective dynamics as well as fast acoustic and gravity modes introduces severe stiffness, strongly constraining the time step size used in fully explicit time integration methods. Efficient numerical discretisations must therefore accurately capture the meteorologically relevant dynamics while reducing the impact of the acoustic stability constraint.

To address the stiffness arising from fast acoustic waves, implicit-explicit (IMEX) time integration methods are particularly attractive for compressible atmospheric flows \cite{giraldo_performance_2024, smolarkiewicz_semi-implicit_2019,melvin_mixed_2019}. In IMEX formulations, the stiff acoustic contributions are treated implicitly, while the slower advective dynamics are advanced explicitly \cite{giraldo_implicit-explicit_2013,benacchio_semi-implicit_2019}, allowing the time step to be determined by the meteorologically relevant flow scales rather than by the speed of sound, while still retaining the fully compressible formulation. For the spatial discretisation, high-order discontinuous Galerkin (DG) methods have been increasingly employed for atmospheric modeling \cite{restelli_semi-lagrangian_2007,giraldo_implicit-explicit_2013,tumolo_semiimplicit_2015,baldauf_comparison_2026}. Their use in atmospheric dynamics is motivated by their ability to combine high-order accuracy with local conservation properties, geometric flexibility, and capability to handle non-conforming adaptive meshes \cite{giraldo_introduction_2020}. Moreover, DG discretisations are particularly well suited to achieving scalability on massively parallel high-performance computing architectures because of their element-local structure \cite{melander_massively_2024,waterhouse_gpu_2026} and efficient matrix-free implementation strategies \cite{kronbichler_fast_2019}.

A high-order IMEX-DG solver for compressible atmospheric flows has been developed over the last years within the \texttt{deal.II} finite element library \cite{orlando_imex-dg_2023,orlando_robust_2024,orlando_impact_2024}. The framework combines IMEX time integration with matrix-free DG operator evaluation, adaptive mesh refinement, and distributed-memory parallelism, enabling scalable simulations of compressible atmospheric dynamics on non-conforming meshes \cite{orlando_efficient_2025,orlando_improving_2025}. However, the current formulation is restricted to non-rotating Cartesian configurations.

For atmospheric applications and numerical weather prediction (NWP), planetary rotation is a leading-order physical effect and a key ingredient of geostrophic balance, i.e. the equilibrium between horizontal pressure-gradient and the Coriolis force that governs much of the large-scale atmospheric circulation \cite{holton_introduction_2013,vallis_atmospheric_2017}. Extending IMEX-DG atmospheric solvers to rotating flows is therefore essential both for standard benchmark configurations involving rotating stratified flows and to proceed towards more realistic large-scale atmospheric applications.

In this work, we extend the existing IMEX-DG framework to include the Coriolis force within the fully compressible Euler equations. The rotational term is incorporated into the momentum equation and consistently included into the implicit part of the IMEX discretisation, together with the pressure gradient and gravity term, while the advective contribution remains explicit.
This generalization modifies the coupled velocity-pressure system arising at each implicit IMEX stage and therefore requires adaptation of the nonlinear solution procedure. In particular, the implicit stage is solved through a nonlinear fixed-point (Picard) iteration \cite{quarteroni_numerical_2007}, in which the pressure subproblem is solved by preconditioned GMRES and the velocity is later recovered through matrix-free operator applications. The introduction of the Coriolis term leads to additional velocity-dependent contributions in the coupled system, whose treatment within the fixed-point iteration is not unique.

We therefore investigate two alternative strategies for incorporating the gravitational and rotational contribution into the nonlinear solution procedure. The two approaches differ in the way the Coriolis term is evaluated during the fixed-point iteration: it can either be lagged using the solution from the previous iteration, or evaluated with the current solution iterate. These alternatives result in different discrete formulations of the implicit stage and may have different implications for convergence of the nonlinear iteration, computational cost and robustness. 
We also analyse the discrete formulation of the rotational contribution and derive the corresponding operator representation for its efficient matrix-free evaluation, ensuring its consistent integration within the IMEX-DG scheme.

The resulting rotating formulation is validated through a sequence of numerical experiments. First, a two-dimensional inertia-gravity wave test in the vertical-slice configuration \cite{baldauf_analytic_2013} is used to verify the implementation and assess its accuracy through convergence studies. These tests are performed on a standard domain and for both rotational treatments. In addition, a planetary-scale inertia-gravity wave simulation is performed to verify the behaviour of the rotating formulation on a much larger spatial scale. 
The study is then extended to fully three-dimensional rotating flows over orography. A bell-shaped mountain configuration is considered to investigate the characteristic effects of rotation in a stratified three-dimensional atmospheric flow. These simulations are carried out both on uniform refined meshes and, to the authors' knowledge for the first time, non-conforming adaptive meshes. Strong scaling tests are also performed to assess the parallel performance of the rotating solver and to compare it with the corresponding non-rotating simulations.
Overall, the numerical results assess the rotating extension in terms of accuracy, robustness and computational efficiency. The results show that the inclusion of Coriolis force preserves the expected accuracy of the IMEX-DG discretisation and enables the solver to reproduce the characteristic effects of rotation in both two- and three-dimensional configurations. The two rotational treatments produce comparable solutions in the tested regimes, while the three-dimensional experiments demonstrate the applicability of the rotating formulation to non-conforming adaptive meshes. The strong scaling results further demonstrate that there is no significant difference in parallel performance between the rotating and non-rotating formulations.

The remainder of the paper is organized as follows. The governing equations are presented in Section~\ref{sec:model}. The temporal and spatial discretisation methods are described in Section~\ref{sec:numerical_framework}. Numerical results in two and three dimensions, including flow over orography and the comparison between the two gravitational and rotational treatments, are presented in Section~\ref{sec:numerical_results}. Some conclusions and considerations about future work are presented in Section \ref{sec:conclusion}.

%% file: Math_model.tex
\section{The mathematical model}
\label{sec:model}

We consider a dry atmosphere modeled as an inviscid fluid governed by the compressible Euler equations on a bounded domain $\Omega \subset \mathbb{R}^d$, with $d=2,3$ \cite{vallis_atmospheric_2017, rieutord_fluid_2015}. The spatial coordinates are denoted by $\bm{x}\in\Omega$ and time by $t \in (0,T_f]$. The system is supplemented by suitable initial and boundary conditions, which will be specified according to the test case considered.

The balance laws of mass, momentum, and total energy read
\begin{subequations}
\begin{align}
\partial_t \rho + \nabla\cdot(\rho\vel) &= 0, \label{eq:euler_mass} \\
\partial_t(\rho\vel)
+ \nabla\cdot(\rho\vel\otimes\vel)
+ \nabla p
&= \rho\bm{a}_{\text{ext}}, \label{eq:euler_mom_rot} \\
\partial_t(\rho E)
+ \nabla\cdot\!\big((\rho E+p)\vel\big)
&= \rho\bm{a}_{\text{ext}}\cdot\vel, \label{eq:euler_energy}
\end{align}
\end{subequations}
where $\rho$ is the density, $\vel=(u,v,w)^{\top}$ is the velocity field, $p$ is the pressure, and $E$ is the specific total energy defined as $E = e + \frac{1}{2}|\vel|^2,$ with $e$ the internal energy per unit mass. Furthermore, $\partial_t$ denotes the time derivative, while $\nabla\cdot\big(\cdots\big)$ and $\nabla\big(\cdots\big)$ denote the three-dimensional divergence and gradient operators, respectively. The system is closed by assuming a calorically perfect ideal gas, for which
\begin{equation}
    p = \rho R T = (\gamma - 1)\rho e,
\end{equation}    
with $\gamma = c_{p}/c_{v}$ the constant ratio of specific heats at constant pressure and constant volume and $R = c_{p} - c_{v}$ the specific gas constant \cite{sandler_chemical_2017}.

We adopt a local Cartesian tangent-plane coordinate system $\bm{x}=(x,y,z)^{\top}$ aligned with the eastward, northward, and upward directions. When expressing the equations in a rotating non-inertial reference frame, additional apparent forces arise \cite{vallis_atmospheric_2017}. Under the $f$-plane approximation, the Coriolis parameter is constant \cite{holton_introduction_2013},
\begin{equation}
f = f_0.
\end{equation}

The centrifugal contribution is absorbed into the definition of the gravitational acceleration, so that the external body force becomes
\begin{equation}
\bm{a}_{\text{ext}} = - g \bm{k} - f\,\bm{k}\times\vel,
\end{equation}
where $\bm{k}$ denotes the unit vector in the vertical direction and $g$ is the gravitational acceleration. The term $-f\,\bm{k}\times\vel$ represents the Coriolis acceleration. Since the Coriolis force is orthogonal to the velocity field, i.e. $(\bm{k} \times \vel)\cdot \vel = 0,$ it does not contribute to the work of external forces in the energy equation.

Using the identity $(\rho E + p)\vel = (h + k)\rho \vel$, where $k = |\vel|^2/2$ and $h$ is the specific enthalpy defined as $h = e + p/\rho$, following the formulation adopted in \cite{orlando_efficient_2022, orlando_imex-dg_2023}, the system can be equivalently rewritten as
\begin{subequations}\label{eq:governing_equations_2}
\begin{align}
\partial_t \rho \;+\; \nabla\!\cdot(\rho \vel) &= 0,
\\[4pt]
\partial_t(\rho \vel)
\;+\; \nabla\!\cdot(\rho \vel\otimes\vel)
\;+\; \nabla p
&= -\,\rho g\,\bm{k} \;-\; \rho f\,\bm{k}\times\vel,
\\[4pt]
\partial_t(\rho E)
\;+\; \nabla\!\cdot\!\big[(h + k)\,\rho \vel\big]
&= -\,\rho g\,\bm{k}\!\cdot\!\vel.
\end{align}
\end{subequations}

For two-dimensional configurations, we employ a vertical slice approximation that restricts the dynamics to the $x-z$ plane, while still allowing for a nonzero meridional velocity component $v$. In practice, this assumption is implemented by imposing  $\partial _y (\cdot)=0$ in \eqref{eq:governing_equations_2}, while the transverse velocity $v$ remains an active prognostic variable. This approach is often used in atmospheric dynamics to reduce the full three-dimensional problem to a two-dimensional slice that still captures the essential balance between advection, pressure gradients, buoyancy, and rotation \cite{benacchio_semi-implicit_2019,melvin_mixed_2019}. Using this approximation, the equations, written in scalar form, become as follows

\begin{subequations}
\begin{align}
\partial_t \rho + \partial_x (\rho u) + \partial_z (\rho w)& = 0,\\
\partial_t (\rho u) + \partial_x (\rho u^2) + \partial_z (\rho u w) 
+  \partial_x p &=f \rho v,\\
\partial_t (\rho v) + \partial_x (\rho u v) + \partial_z (\rho v w) 
 &= - f \rho u,\\
\partial_t (\rho w) + \partial_x (\rho u w) + \partial_z (\rho w^2) 
+  \partial_z p &= -g\rho,\\
\partial_t (\rho E) 
+ \partial_x \big[(h+k)\rho u \big] 
+ \partial_z \big[(h+k)\rho w \big] 
&= - g \rho w.
\label{eq:vertical_slice}
\end{align}
\end{subequations}

In the present work, following \cite{cotter_compatible_2023,schmid_numerical_2021}, the Coriolis term is formulated in terms of the velocity perturbation $\vel'=\vel-\vel_0$, where $\vel_0$ denotes the prescribed background velocity. This choice removes the Coriolis contribution associated with the background flow, which, in geostrophic balance, should be compensated by a transverse pressure gradient that cannot be represented in a vertical slice configuration. Accordingly, the momentum equation in \eqref{eq:governing_equations_2} is modified as 
\begin{equation}
    \partial_t(\rho\vel)
+\nabla\cdot(\rho\vel\otimes\vel)
+\nabla p
=-\rho g\bm{k}-\rho f\bm{k}\times\vel'.
\end{equation}
The same formulation is retained in the three-dimensional simulations (Section \ref{sec:3D}), where periodic boundary conditions in the $y$-direction may likewise prevent the representation of the pressure gradient required to balance the Coriolis force associated with the background flow.

%% file: Numerical_Framework.tex
\section{The numerical framework}
\label{sec:numerical_framework}

The numerical approximation of the governing equations \eqref{eq:governing_equations_2} is based on a high-order DG spatial discretization coupled with an IMEX time integration strategy. This combination is well suited for atmospheric flows, where slow advective dynamics coexist with fast acoustic waves. In the low Mach number regime, terms proportional to the inverse of the squared Mach number in the non-dimensional formulation of the equations, such as the pressure gradient, introduce stiffness, leading to severe stability restrictions for fully explicit schemes \cite{ascher_implicit-explicit_1997,giraldo_study_2008,gassner_novel_2021}.

Following \cite{orlando_efficient_2022,orlando_imex-dg_2023}, time discretization is performed using an IMEX additive Runge--Kutta (ARK) method, in which stiff contributions are treated implicitly and non-stiff terms explicitly \cite{hairer_solving_1996, boscarino_implicit-explicit_2025, pareschi_implicitexplicit_2005}. In particular, advective and nonlinear transport terms are handled explicitly, while pressure-gradient terms, together with gravity and rotation, are included in the implicit part. Although the Coriolis term and the gravity term are not intrinsically stiff, they are commonly included in the implicit treatment to better preserve the hydrostatic and geostrophic balances.

In this work, we employ a second-order IMEX-ARK scheme, defined by two sets of coefficients, associated with the explicit and implicit parts, respectively. The corresponding Butcher tableaux are reported in Table \ref{tab:butcher_tableaux}.

\begin{table}[h!]
\centering
\begin{minipage}{0.45\textwidth}
\centering
\[
\begin{array}{c|ccc}
0 & 0 &  &  \\
\chi & \chi & 0 &  \\
1 & 1-a_{32} & a_{32} & 0 \\
\hline
& \frac{1}{2}-\frac{\chi}{4}
& \frac{1}{2}-\frac{\chi}{4}
& \frac{\chi}{2}
\end{array}
\]
\end{minipage}
\hfill
\begin{minipage}{0.45\textwidth}
\centering
\[
\begin{array}{c|ccc}
0 & 0 &  &  \\
\chi & \frac{\chi}{2} & \frac{\chi}{2} &  \\
1 & \frac{1}{2}-\frac{\chi}{4}
& \frac{1}{2}-\frac{\chi}{4}
& \frac{\chi}{2} \\
\hline
& \frac{1}{2}-\frac{\chi}{4}
& \frac{1}{2}-\frac{\chi}{4}
& \frac{\chi}{2}
\end{array}
\]
\end{minipage}
\caption{Butcher tableaux for the explicit (left) and implicit (right) RK methods.}
\label{tab:butcher_tableaux}
\end{table}

For the implicit part, we adopt the TR-BDF2 scheme as the implicit companion method \cite{hosea_analysis_1996}. The parameter $\chi$ is set to $\chi=2-\sqrt{2}$, which ensures $L$-stability and thus provides effective damping of fast acoustic modes in stiff low-Mach regimes \cite{hosea_analysis_1996,tumolo_semiimplicit_2015,bonaventura_unconditionally_2017}. The implicit tableau has an explicit first stage, followed by singly diagonally implicit stages with identical diagonal coefficients. It therefore has the structure of an ESDIRK method \cite{kennedy_higher-order_2019, ketcheson_optimal_2009}.

For the explicit part, we employ a second-order RK scheme \cite{giraldo_implicit-explicit_2013}. The coefficients are chosen such that $b=\tilde{b}$ and $c=\tilde{c}$,  where $b,\widetilde{b}$ denote the final update weights and $c,\widetilde{c}$ denote the stage nodes of the explicit and implicit Runge-Kutta tableaux, respectively \cite{boscarino_implicit-explicit_2025}. The condition $b=\tilde{b}$ guarantees preservation of linear invariants \cite{giraldo_implicit-explicit_2013}, while $c=\tilde{c}$ is usually employed to simplify the order conditions. The remaining free parameter in the explicit tableau is set to $a_{32}=0.5$, following \cite{orlando_efficient_2022}, so as to improve monotonicity properties.

Thus, for an ODE of the form $\mathbf{y}'=\mathbf{f}_{\mathrm{S}}(\mathbf{y},t)+\mathbf{f}_{\mathrm{NS}}(\mathbf{y},t)$, where NS denotes the non-stiff term and S denotes the stiff one, the $l$-th stage of the solution is computed as
\begin{equation}
\begin{aligned}
\mathbf{y}^{(n,l)} = \mathbf{y}^{n}
&+ \Delta t \sum_{m=1}^{l-1}
\Big[
a_{lm}\mathbf{f}_{\mathrm{NS}}(\mathbf{y}^{(n,m)},t^n+c_m\Delta t) 
+ \widetilde{a}_{lm}\mathbf{f}_{\mathrm{S}}(\mathbf{y}^{(n,m)},t^n+c_m\Delta t)
\Big] \\&+ \Delta t \, \widetilde{a}_{ll} \,
\mathbf{f}_{\mathrm{S}}(\mathbf{y}^{(n,l)},t^n+c_l\Delta t).
\end{aligned}
\end{equation}
for $l=1,2,3$, where $a_{lm}$ and  $\tilde{a}_{lm}$ denote the coefficients of the explicit and implicit Butcher tableaux, respectively.
The update is then given by
\begin{equation}
\mathbf{y}^{n+1} = \mathbf{y}^n
+ \Delta t \sum_{l=1}^{3} b_l \left[
\mathbf{f}_{\mathrm{NS}}(\mathbf{y}^{(n,l)}, t^n + c_l \Delta t)
+ \mathbf{f}_{\mathrm{S}}(\mathbf{y}^{(n,l)}, t^n + c_l \Delta t)
\right].
\end{equation}

Applying the IMEX-ARK method presented in Table \ref{tab:butcher_tableaux}, the first stage is purely formal and  simply reduces to $\rho^{(n,1)}=\rho^n$ for the density and analogously for the other prognostic variables. 
From the second stage onward, it is possible to write the following general formulae
\begin{subequations}{
\begin{align}
\rho^{(n,l)} &= \rho^n - \sum_{m=1}^{l-1}a_{lm}\Delta t \, \nabla \cdot (\rho^{(n,m)} \vel^{(n,m)}),  \\
\rho^{(n,l)} \vel^{(n,l)} & + \tilde{a}_{ll}\Delta t \nabla p^{(n,l)} +g \tilde{a}_{ll}\Delta t \rho^{(n,l)} \bm{k}+ f\tilde{a}_{ll} \Delta t \rho^{(n,l)}\bm{k}\times \vel^{(n,l)}  = \mathbf{m}^{(n,l)}, \\
\rho^{(n,l)} E^{(n,l)} &+ \tilde{a}_{ll}\Delta t \, \nabla \cdot (h^{(n,l)} \rho^{(n,l)} \vel^{(n,l)}) + g\,\tilde{a}_{ll}\rho^{(n,l)} \bm{k} \cdot \vel^{(n,l)}=\hat{e}^{(n,l)},
\end{align}
\label{eq:IMEXsystem_sstage}
}
\end{subequations}
where
\begin{equation}
\begin{aligned}
\mathbf{m}^{(n,l)} = \rho^n \vel^n - \sum_{m=1}^{l-1} [&a_{lm}\Delta t \, \nabla \cdot (\rho^{(n,m)} \vel^{(n,m)} \otimes \vel^{(n,m)})+ \tilde{a}_{lm}\Delta t \nabla p^{(n,m)}  + g \tilde{a}_{lm}\Delta t \rho^{(n,m)} \bm{k} \\&\quad+f\tilde{a}_{lm}\Delta t \rho^{(n,m)} \bm{k}\times (\vel^{(n,m)}-\vel_0)]+f \tilde a_{ll} \Delta t \rho^{(n,l)} \bm{k} \times \vel_0,
\end{aligned}
\end{equation}
\begin{equation}
\begin{aligned}
\hat{e}^{(n,l)} = \rho^n E^n - \sum_{m=1}^{l-1} &\Bigg[\tilde{a}_{lm}\Delta t \, \nabla \cdot (h^{(n,m)} \rho^{(n,m)} \vel^{(n,m)}) + a_{lm}\Delta t  \,\nabla \cdot (\kappa^{(n,m)} \rho^{(n,m)} \vel^{(n,m)})  \\
 &+g\Delta t \,\tilde{a}_{lm}\rho^{(n,m)} \bm{k} \cdot \vel^{(n,m)}\Bigg]. 
\end{aligned}
\end{equation}
\\

The spatial discretization of the governing equations is based on a high-order nodal DG method \cite{hesthaven_nodal_2008} on a mesh $\mathcal{T}_h$, composed of quadrilateral elements in 2D and hexaedral elements in 3D. This choice combines high-order accuracy with local conservation and geometric flexibility, making it well suited for compressible atmospheric flows \cite{hesthaven_nodal_2008}.

For an integer $r \geq 0$, we define the polynomial space $\mathbb{Q}_r(K)$ as the space of polynomials of degree $\leq r$ in each coordinate direction on the element $K$.
The global DG finite element space is then
\[
Q_r = \{ v \in L^2(\Omega) \ : \ v|_K \in \mathbb{Q}_r(K), \ \forall K \in \mathcal{T}_h \}.
\]
Scalar quantities such as density and pressure are approximated in $Q_r$, 
while vector-valued quantities such as velocity in $[Q_r]^d= \{ v = (v_1,\dots,v_d): v_i \in Q_r \},$ where we recall that $d$ denotes the number of the physical dimensions. We employ a nodal DG formulation based on tensor-product elements. On each element $K$, the discrete solution is
\[
\mathbf{q}_h\big|_K(\bm{x},t) = \sum_{j=1}^{N_p} \mathbf{q}_j(t) \ell_j(\bm{x}),
\]
where $\{\mathbf{q}_j \approx \mathbf{q}(\bm{x}_j)\}$, $j=1,\dots,N_{p}$, with $N_p=(r+1)^d$, are nodal values at Gauss-Legendre--Lobatto points $\{\bm{x}_j\}$, and $\{\ell_j\}$ are the associated Lagrange basis functions~\cite{quarteroni_numerical_2007}. The tensor-product structure enables efficient sum-factorized quadrature and operator evaluation on quadrilateral and hexahedral meshes~\cite{hesthaven_nodal_2008,kronbichler_fast_2019}.

Following the standard DG formulation, the equations are multiplied by test functions and integrated by parts element-wise, leading to surface flux contributions that couple neighbouring elements. A Rusanov-type flux is used at interfaces, with stabilization parameter $\lambda = |\mathbf{u}\cdot\mathbf{n}|$ \cite{orlando_efficient_2022}. For the pressure-gradient contribution, centered fluxes are employed in order to simplify the implicit treatment within the IMEX formulation, as implicit time integration inherently ensures stability without the need for upwind stabilization.

Denoting $\{\boldsymbol{\varphi}_i\}_{i=1}^{\dim([Q_r]^d)}$ the basis functions of the DG velocity space $[Q_r]^d$ and $\{\psi_j\}_{j=1}^{\dim(Q_r)}$ the basis functions of the DG pressure space ${Q}_r$, velocity and pressure can be rewritten in a discrete form as
 \begin{equation}
 \begin{aligned}
     \vel(\bm{x},t)\approx\sum_{j=1}^{\text{dim}([{Q}_r]^d)} u_j(t) \boldsymbol{\varphi}_j(\bm{x}), & \qquad p(\bm{x},t)\approx \sum_{j=1}^{\text{dim}(Q_r)} p_j(t) {\psi}_j(\bm{x}).
     \end{aligned}
 \end{equation}
 
Given these definitions, it is possible to derive the weak formulation of the momentum equation at stage $l$ in matrix form:  

\begin{equation}
\mathbf{A}^{(n,l)}\,\mathbf{u}^{(n,l)} + \mathbf{B}^{(n,l)} \mathbf{p}^{(n,l)} + \mathbf{R}^{(n,l)}\mathbf{u}^{(n,l)}= \mathbf{f}^{(n,l)}- \mathbf{f}_g^{(n,l)},
\label{eq:momentum_matrix}
\end{equation}
where
$\mathbf{u}^{(n,l)}$ and $\mathbf{p}^{(n,l)}$ are the vectors of velocity and pressure degrees of freedom,

\begin{subequations}
\begin{equation}
\begin{aligned}
A^{(n,l)}_{ij}
=&
\sum_K \int_K
\rho^{(n,l)}
\boldsymbol{\varphi}_j \cdot \boldsymbol{\varphi}_i
\, \mathrm{d}\Omega,
\end{aligned}
\end{equation}
\begin{equation}
\begin{aligned}
B^{(n,l)}_{ij}
=&
\sum_K \int_K
-\tilde{a}_{ll}\Delta t\,
\nabla \cdot \boldsymbol{\varphi}_i\,\psi_j
\,\mathrm{d}\Omega+
\sum_{\Gamma \in E} \int_\Gamma
\tilde{a}_{ll}\Delta t\,
\avg{\psi_j}
\jump{\boldsymbol{\varphi}_i}
\,\mathrm{d}\Sigma,
\end{aligned}
\end{equation}
\begin{equation}
\begin{aligned}
 R^{(n,l)}_{ij}=&  \tilde{a}_{ll}\Delta t f
\sum_{K\in\mathcal{T}_h}\!\int_K 
 \rho^{(n,l)}
 (
\bm{k} \times 
 \boldsymbol{\varphi}_j
 )
\cdot \boldsymbol{\varphi}_i\, \mathrm{d} \Omega, \quad \qquad\\
\mathbf{f}_{g_{i}}^{(n,l)}=&\,
\tilde{a}_{ll}\,\Delta t\,g
\sum_{K\in\mathcal{T}_h}\!\int_K 
  \rho^{(n,l)}\bm{k}
\cdot \boldsymbol{\varphi}_{i}\, \mathrm{d} \Omega,
\end{aligned}
\end{equation}
\begin{equation}
\begin{aligned}
\mathbf{f}^{(n,l)} =&
\sum_{K \in \mathcal{T}_h} \int_K \rho^n \vel^n \cdot \boldsymbol{\varphi}_i \, \mathrm{d} \Omega
+ \sum_{m=1}^{l-1}\Bigg[ \sum_{K \in \mathcal{T}_h} \int_K a_{lm}\Delta t \, \left(\rho^{(n,m)} \vel^{(n,m)} \otimes \vel^{(n,m)}\right) : \nabla \boldsymbol{\varphi}_i \, \mathrm{d}\Omega 
\\
+& \sum_{K \in \mathcal{T}_h} \int_K \tilde{a}_{lm}\Delta t p^{(n,m)} \nabla\!\cdot \boldsymbol{\varphi}_i  \mathrm{d} \Omega
- \sum_{\Gamma \in \mathcal{E}} \int_\Gamma  a_{lm}\Delta t \, \avg{ \rho^{(n,m)} \vel^{(n,m)} \otimes \vel^{(n,m)}} : \tensorjump{\boldsymbol{\varphi}_i} \, \mathrm{d} \Sigma \nonumber
\\
-& \sum_{\Gamma \in \mathcal{E}} \int_\Gamma  \tilde{a}_{lm}\Delta t \, \avg{ p^{(n,m)}} \, \jump{\boldsymbol{\varphi}_i} \, \mathrm{d} \Sigma
-\sum_{\Gamma \in \mathcal{E}} \int_\Gamma a_{lm}\Delta t \, \frac{\lambda^{(n,m)}}{2} \, \jump{ \rho^{(n,m)} \vel^{(n,m)} } : \tensorjump{\boldsymbol{\varphi}_i} \, \mathrm{d} \Sigma\nonumber
\\ 
- &\sum_{K \in \mathcal{T}_h} \int_K g \tilde{a}_{lm}\Delta t \, \rho^{(n,m)}  \bm{k} \cdot \boldsymbol{\varphi}_i \, \mathrm{d} \Omega
 - \sum_{K \in \mathcal{T}_h}
\int_K \tilde{a}_{lm}\Delta t f \rho^{(n,m)}(\bm{k}\times (\vel^{(n,m)}-\vel_0))\cdot\boldsymbol{\varphi}_i \mathrm{d} \Omega \Bigg]\nonumber\\+& \sum_{K \in \mathcal{T}_h} \int_K \tilde a_{ll} \Delta t f\rho^{(n,l)}(\bm{k} \times  \boldsymbol{\vel}_{0})\cdot \boldsymbol{\varphi}_i\, \mathrm{d} \Omega, \nonumber
\end{aligned}
\end{equation}
\end{subequations}
with
$\lambda^{(n,m)} = \max \left( | \vel^{(n,m),+}\cdot \mathbf{n}^+|, \; | \vel^{(n,m),-}\cdot \mathbf{n}^- | \right),$ and the symbol $:$ indicates the scalar product for second-order tensors.

Taking into account that
\begin{equation}
\rho^{(n,l)} E^{(n,l)} = \rho^{(n,l)} e^{(n,l)}(p^{(n,l)},\rho^{(n,l)}) + \rho^{(n,l)} \kappa^{(n,l)},
\end{equation}
it is also possible to derive the algebraic form of the energy equation as

\begin{equation}
\mathbf{C}^{(n,l)} \mathbf{u}^{(n,l)} + \mathbf{D}^{(n,l)} \mathbf{p}^{(n,l)}+\mathbf{k}^{(n,l)}+\mathbf{M}^{(n,l)}_g \mathbf{u}^{(n,l)}= \mathbf{g}^{(n,l)},
\end{equation}
where 
\begin{subequations}
\begin{equation}
\begin{aligned}
C^{(n,l)}_{ij} = &\sum_{K\in \mathcal{T}_h} \int_K -\tilde{a}_{ll}\Delta t \, h^{(n,l)} \rho^{(n,l)} \boldsymbol{\varphi}_j \cdot \nabla \psi_i \, \mathrm{d} \Omega
+ \sum_{\Gamma \in E} \int_\Gamma \tilde{a}_{ll}\Delta t \, \{\!\! \{ h^{(n,l)} \rho^{(n,l)} \boldsymbol{\varphi}_j \}\!\! \} \cdot \jump{\psi_i} \, \mathrm{d} \Sigma,
\end{aligned}
\end{equation}
\begin{equation}
    D^{(n,l)}_{ij}=\sum_{K\in \mathcal{T}_h}\int_K \frac{1}{(\gamma -1)}\psi_i\psi_j \mathrm{d}\Omega ,
\end{equation}

\begin{equation}
    {k}^{(n,l)}_i=\sum_{K \in \mathcal{T}_h}\int_K \rho^{(n,l)} k^{(n,l)}\mathbf{\psi}_i\mathrm{d}\Omega,
\end{equation}

\begin{equation}
[{M}_{g}^{(n,l)}]_{ij}=\sum_{K\in\mathcal{T}_h}\int_K g\tilde{a}_{ll}\Delta t\,\rho^{(n,l)}
\bm{k}\cdot\varphi_i\,\psi_j\,\mathrm{d}\Omega
\end{equation}
and the expression of $\mathbf{g}$ is given by

\begin{equation}
\begin{aligned}
\mathbf{g}^{(n,l)}
&= \sum_{K\in\mathcal{T}_h}\int_K \rho^n E^n \psi_i \,\mathrm{d}\Omega
+\sum_{m=1}^{l-1}\Bigg[
\sum_{K\in\mathcal{T}_h}\int_K \Delta t
a_{lm}\,\left(\kappa^{(n,m)} \rho^{(n,m)}\vel^{(n,m)}\right)\cdot\nabla\psi_i\,\mathrm{d}\Omega\\
& \quad +\sum_{K\in\mathcal{T}_h}\int_K \Delta t\tilde{a}_{lm}\left(h^{(n,m)}\rho^{(n,m)}\vel^{(n,m)}\right)
\cdot\nabla\psi_i\,\mathrm{d}\Omega
\\
&\quad
-\sum_{\Gamma\in\mathcal{E}}\int_\Gamma \Delta t\Big(
a_{lm}\,\avg{\kappa^{(n,m)}\rho^{(n,m)}\vel^{(n,m)}}
+\tilde{a}_{lm}\avg{h^{(n,m)}\rho^{(n,m)}\vel^{(n,m)}}
\Big)\cdot\jump{\psi_i}\,\mathrm{d}\Sigma
\\
&\quad
-\sum_{\Gamma\in\mathcal{E}}\int_\Gamma 
\Delta t\,\frac{\lambda^{(n,m)}}{2}
\Big(a_{lm}\jump{\rho^{(n,m)}\kappa^{(n,m)}}
+\tilde{a}_{lm}\jump{\rho^{(n,m)}e^{(n,m)}}\Big)\cdot\jump{\psi_i}\,\mathrm{d}\Sigma
\\
&\quad
-\sum_{K\in\mathcal{T}_h}\int_K g
\,\tilde{a}_{lm}\Delta t\,\rho^{(n,m)}\bm{k}\cdot\vel^{(n,m)}\psi_i \,\mathrm{d}\Omega\\ 
&\quad- \frac{\lambda^{(n,m)}\tilde a_{lm}\Delta t}{2}\sum_{\Gamma \in \mathcal{E}}\int_\Gamma \jump{\rho^{(n,m)}E^{(n,m)}}\cdot\jump{\psi_i}\mathrm{d}\Sigma \Bigg]\\&\quad-\frac{\lambda^{(n,l)}\tilde a_{ll}\Delta t}{2}\sum_{\Gamma \in \mathcal{E}}\int_\Gamma\jump{\rho^{(n,l)}e^{(n,l)}\psi_i}\cdot\jump{\psi_j}\mathrm{d}\Sigma.
\end{aligned}
\end{equation}
with
\[
\lambda^{(n,m)} = \max \left( | \vel^{(n,m),+}\cdot \mathbf{n}^+|, \; | \vel^{(n,m),-}\cdot \mathbf{n}^- | \right).
\]
\end{subequations}
Then, we end up with the following coupled system of equations
\begin{equation}
\begin{aligned}
&\mathbf{A}^{(n,l)}\,\mathbf{u}^{(n,l)}+\mathbf{R}^{(n,l)}\mathbf{u}^{(n,l)} + \mathbf{B}^{(n,l)} \mathbf{p}^{(n,l)} + \mathbf{f}_g^{(n,l)}= \mathbf{f}^{(n,l)},\\
&\mathbf{C}^{(n,l)} \mathbf{u}^{(n,l)} +\mathbf{M}^{(n,l)}_g \mathbf{u}^{(n,l)} + \mathbf{D}^{(n,l)} \mathbf{p}^{(n,l)}+\mathbf{k}^{(n,l)}= \mathbf{g}^{(n,l)}.
\label{eq:matrix_form}
\end{aligned}
\end{equation}

The algebraic system for $(\mathbf{u}^{(n,l)},\mathbf{p}^{(n,l)})$ is nonlinear and coupled, as matrices and right-hand sides depend on the stage unknowns through thermodynamic quantities. Following the strategy commonly adopted in semi-implicit atmospheric solvers \cite{casulli_pressure_1984,orlando_efficient_2022,orlando_imex-dg_2023}, the nonlinear system is solved through a fixed-point (Picard) iteration \cite{quarteroni_numerical_2007,dumbser_conservative_2016}. This approach avoids the solution of a fully nonlinear system at each IMEX stage while preserving the implicit coupling between pressure and velocity. At each iteration, selected terms are evaluated using the previous iterate, resulting in a sequence of linear systems. The iteration is initialized using the solution obtained at the previous IMEX stage and proceeds until the relative variation between two consecutive iterates falls below a prescribed tolerance. Different lagging strategies for the nonlinear terms result in separate formulations:

\begin{itemize}
\item In the first approach (R1), the term $\mathbf{R}^{(n,l)}\mathbf{u}^{(n,l)}$ in the momentum equation is evaluated at the previous Picard iterate $k$, while the remaining terms are treated at the current iteration $k+1$. The same strategy is adopted for the gravitational contribution in the energy equation, where $\mathbf{M}_g^{(n,l)}\mathbf{u}^{(n,l)}$ and $\mathbf{k}^{(n,l)}$ are lagged at iteration $k$. This choice moves the rotational and gravitational contributions to the right-hand side, where they are evaluated using the previous Picard iterate. The resulting system reads
\begin{equation}{
\begin{aligned}
&\mathbf{A}^{(n,l)}\,\mathbf{u}^{(n,l,k+1)}+ \mathbf{B}^{(n,l)} \mathbf{p}^{(n,l,k+1)} = \mathbf{f}^{(n,l)}-  \mathbf{f}_g^{(n,l)}- \mathbf{R}^{(n,l)}\mathbf{u}^{(n,l,k)} ,\\
&\mathbf{C}^{(n,l,k)}\mathbf{u}^{(n,l,k+1)}+\mathbf{D}^{(n,l)}\mathbf{p}^{(n,l,k+1)}=\mathbf{g}^{(n,l,k)}-\mathbf{k}^{(n,l,k)}-\mathbf{M}_g^{(n,l)}\mathbf{u}^{(n,l,k)}.
\label{eq:matrix_form_A_C}
\end{aligned}
}
\end{equation}

This coupled system is solved using a Schur complement approach. Eliminating the velocity from the momentum equation yields a reduced system for the pressure,
\begin{equation}
\left[\mathbf{D}^{(n,l)}
- \mathbf{C}^{(n,l,k)} \big(\mathbf{A}^{(n,l)}\big)^{-1}
  \mathbf{B}^{(n,l)}\right]\,
\mathbf{p}^{(n,l,k+1)}
=
\mathbf{RHS}^{(n,l,k)}_{R1},
\label{eq:matrix_form_R1}
\end{equation}
where
\begin{equation}
\begin{aligned}
\mathbf{RHS}^{(n,l,k)}_{R1}
=&\;\mathbf{g}^{(n,l,k)}
- \mathbf{M}_g^{(n,l)}\mathbf{u}^{(n,l,k)}
- \mathbf{k}^{(n,l,k)} \\&\quad- \mathbf{C}^{(n,l)} \big(\mathbf{A}^{(n,l)}\big)^{-1}
\Big(\mathbf{f}^{(n,l)}
- \mathbf{R}^{(n,l)}\mathbf{u}^{(n,l,k)}
- \mathbf{f}_g^{(n,l)}\Big).
\end{aligned}
\end{equation}

\item In the second strategy (R2), both the rotational term $\mathbf{R}^{(n,l)}\mathbf{u}^{(n,l)}$ and the gravitational contribution $\mathbf{M}_g^{(n,l)}\mathbf{u}^{(n,l)}$ are evaluated at the current Picard iteration $k+1$ and retained on the left-hand side. These contributions are therefore treated implicitly within the coupled momentum-energy system.

Proceeding as before, the resulting linear system for the pressure reads
\begin{equation}
\big(\mathbf{D}^{(n,l)}
- (\mathbf{C}^{(n,l,k)}+\mathbf{M}_g^{(n,l)}) \big(\mathbf{A}^{(n,l)}+\mathbf{R}^{(n,l)}\big)^{-1}
  \mathbf{B}^{(n,l)}\big)\,
\mathbf{p}^{(n,l,k+1)}
=
\mathbf{RHS}^{(n,l,k)}_{R2},
\label{eq:matrix_form_R2}
\end{equation}
with
\begin{equation}
\begin{aligned}
\mathbf{RHS}^{(n,l,k)}_{R2}
=&\;\mathbf{g}^{(n,l,k)}
- \mathbf{k}^{(n,l,k)}  - (\mathbf{C}^{(n,l,k)}+\mathbf{M}_g^{(n,l)}) \big(\mathbf{A}^{(n,l)}+\mathbf{R}^{(n,l)}\big)^{-1}
\Big(\mathbf{f}^{(n,l)}
- \mathbf{f}_g^{(n,l)}\Big).
\end{aligned}
\end{equation}
Once $\mathbf{p}^{(n,l,k+1)}$ has been computed, the velocity $\mathbf{u}^{(n,l,k+1)}$ is recovered from the momentum equation.
\end{itemize}

In both approaches, the pressure system is solved using a GMRES iterative method, with a preconditioning strategy described in \cite{orlando_improving_2025}.

\subsection{Matrix-free application of $(\mathbf{A}+\mathbf{R})^{-1}$}
\label{sec:matrix_free}

The formulations discussed above require repeated applications of block-diagonal operators such as $\mathbf{A}^{-1}$ in \ref{eq:matrix_form_R1} and $(\mathbf{A}+\mathbf{R})^{-1}$ in \ref{eq:matrix_form_R2}. These operators are evaluated in a matrix-free fashion by exploiting the tensor-product structure of the DG discretisation on quadrilateral/hexahedral meshes \cite{kronbichler_fast_2019,arpaia_highorder_2026}.

The density-weighted mass matrix can be expressed element-wise as
\begin{equation}
\mathbf{A}_K = \mathbf{S}^T \mathbf{W}_\rho \mathbf{S},
\end{equation}
where $\mathbf{S}$ is the interpolation matrix between the nodal and quadrature bases, and $\mathbf{W}_\rho$ is a diagonal matrix containing quadrature weights, geometric factors, and density \cite{kronbichler_fast_2019}. This factorisation enables an efficient application of $\mathbf{A}_K^{-1}$ through tensor-product operations.

For the R2 strategy, the inverse of the modified operator $(\mathbf{A}+\mathbf{R})$ is required. Rather than treating $\mathbf{R}$ as a generic matrix contribution, we exploit the algebraic structure of the Coriolis operator. Introducing the linear mapping $\mathcal{J}\mathbf{v} = \mathbf{k}\times\mathbf{v}$, one can show that
\begin{subequations}
 
\begin{equation}
\mathbf{R} = \beta\, \mathbf{A}\mathbf{J},
\qquad
\beta = \tilde a_{ll}\Delta t f,
\end{equation}
so that
\begin{equation}
\mathbf{A}+\mathbf{R} = \mathbf{A}(\mathbf{I}+\beta\mathbf{J}),
\end{equation}
under the assumption that $f$ is constant.
This immediately yields the factorisation
\begin{equation}
(\mathbf{A}+\mathbf{R})^{-1} = (\mathbf{I}+\beta\mathbf{J})^{-1}\mathbf{A}^{-1}.
\end{equation}
\end{subequations}
Assuming a constant vertical rotation axis $\mathbf{k}=(0,0,1)^T$, the operator $\mathcal{J}$ acts as a planar rotation on the horizontal velocity components, satisfying $\mathcal{J}^2=-I$. Hence,
\begin{equation}
(\mathbf{I}+\beta\mathbf{J})^{-1}
=
\frac{1}{1+\beta^2}
\begin{pmatrix}
1 & \beta & 0 \\
-\beta & 1 & 0 \\
0 & 0 & 1+\beta^2
\end{pmatrix},
\end{equation}
which corresponds to a local transformation applied pointwise to the velocity field.
Therefore, the application of $(\mathbf{A}+\mathbf{R})^{-1}$ reduces to a standard matrix-free application of $\mathbf{A}^{-1}$, followed by a local algebraic transformation involving $(\mathbf{I}+\beta\mathbf{J})^{-1}$.
This avoids the need to solve linear systems with the full coupled operator and preserves the efficiency of the matrix-free implementation.

It is worth noting that the factorisation introduced above is only valid if $f$ is constant. If $f$ varies spatially, the more general identity $(\mathbf{A}+\mathbf{R})^{-1}
=(\mathbf{I}+\mathbf{A}^{-1}\mathbf{R})^{-1}\mathbf{A}^{-1}$ remains valid. Since $\mathbf{A}$ and $\mathbf{R}$ are both block-diagonal with respect to the DG elements, $\mathbf{A}^{-1}\mathbf{R}$ retains this structure, allowing the inverse to be evaluated locally without assembling a global matrix. The detailed treatment of spatially varying rotation parameters is beyond the scope of the present work and is left for future investigation.

%% file: Numerical_results.tex
\section{Numerical results}
\label{sec:numerical_results}

The numerical method presented in the previous sections is validated through a series of two- and three-dimensional test cases for rotating atmospheric flows. The considered benchmarks are designed to evaluate the ability of the proposed IMEX-DG formulation to accurately capture wave propagation and balanced rotating dynamics while preserving long-time stability. The numerical experiments include both standard rotating inertia-gravity wave benchmarks and fully three-dimensional simulations of stratified flow over orography. The latter provide a more challenging setting to assess the behaviour of the method in realistic large-scale atmospheric regimes. Simulations are performed on both uniform and non-conforming meshes in order to evaluate the robustness of the formulation in the presence of local mesh refinement.

Particular attention is devoted to verifying that the inclusion of rotational effects preserves the accuracy, robustness, and scalability properties of the original non-rotating IMEX-DG solver implemented within the matrix-free distributed-memory framework of the \texttt{deal.II} finite element library \cite{bangerth_dealiigeneral-purpose_2007}. 

For the numerical experiments, we define the minimum mesh element diameter $\mathcal{H}$ as 
\begin{equation}
    \mathcal{H}=\min{\{\text{diam}(K)|K \in\mathcal{T}_h}\},
\end{equation}
as well as the acoustic Courant number $C$ and advective Courant number $C_{\mathrm{adv}}$ as
\begin{equation}
    C=\frac{rc \Delta t}{\mathcal{H}}\sqrt{d}, \qquad C_{\mathrm{adv}}=\frac{r U \Delta t}{\mathcal{H}}\sqrt{d},
\end{equation}
where $c$ is the reference magnitude of speed of sound, $U$ is the reference magnitude of the flow velocity, $r$ is  the polynomial degree of the DG method and  $\mathcal{H}$ is interpreted as an effective spatial resolution for a DG discretisation \cite{orlando_imex-dg_2023}.  
In all the simulations, we set $\gamma=1.4$, $g=9.81\text{m s}^{-2}$, and $R=287 \text{  J Kg}^{-1}\text{K}^{-1}$.

Simulations have been performed on computing clusters at the Technical University of Denmark on 28 Lenovo ThinkSystem SD530 compute nodes equipped with two Intel Xeon Gold 6226R processors (16 physical cores, 2.90 GHz), \cite{noauthor_dtu_2026} and on the MeluXina HPCF, using up to 1024 2x AMD EPYC Rome 7H12 64c 2.6GHz CPUs \cite{noauthor_meluxina_2026}, compiled with GCC~12.3, targeting 256-bit vectorization.

\subsection{Two-dimensional inertia-gravity wave}
\begin{figure}[h!]
  \centering

  \begin{subfigure}{0.65\textwidth}

    \centering
    \includegraphics[width=\linewidth]{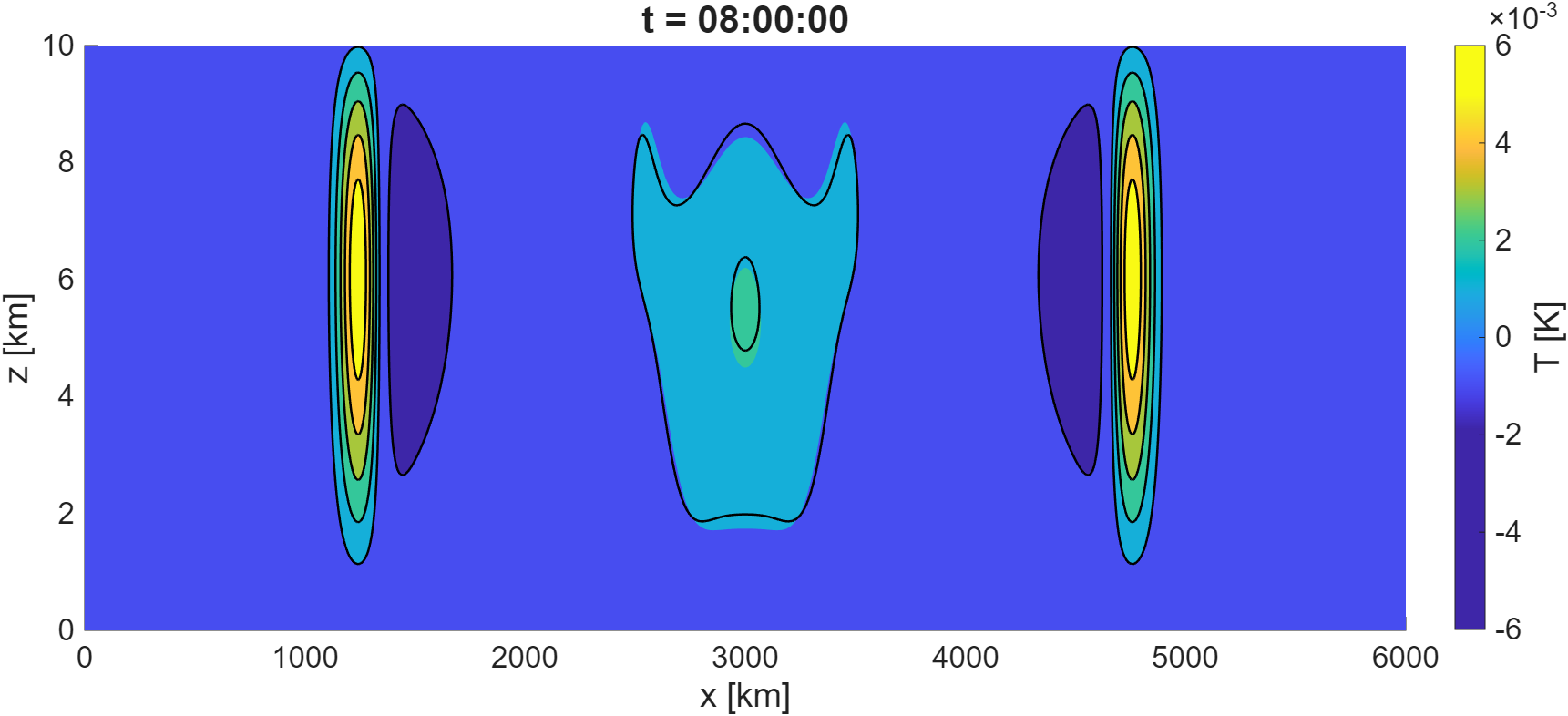}
    \caption{}
    \label{fig:Baldauf_T}
  \end{subfigure}

  \vspace{0.2cm}

  \begin{subfigure}{0.65\textwidth}
    \centering
    \includegraphics[width=\linewidth]{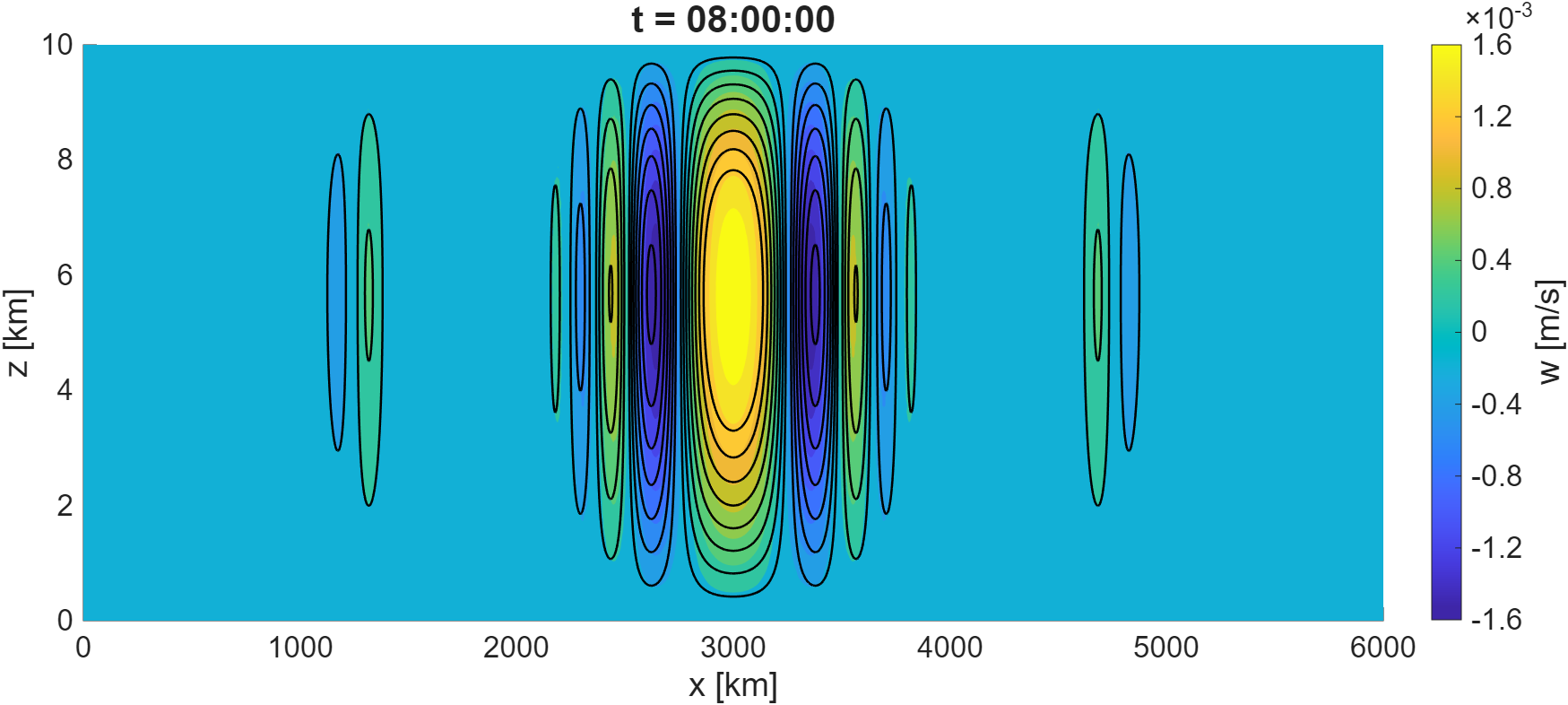}
    \caption{}
    \label{fig:Baldauf_w}
  \end{subfigure}

  \vspace{0.2cm}

  \begin{subfigure}{0.65\textwidth}
    \centering
    \includegraphics[width=\linewidth]{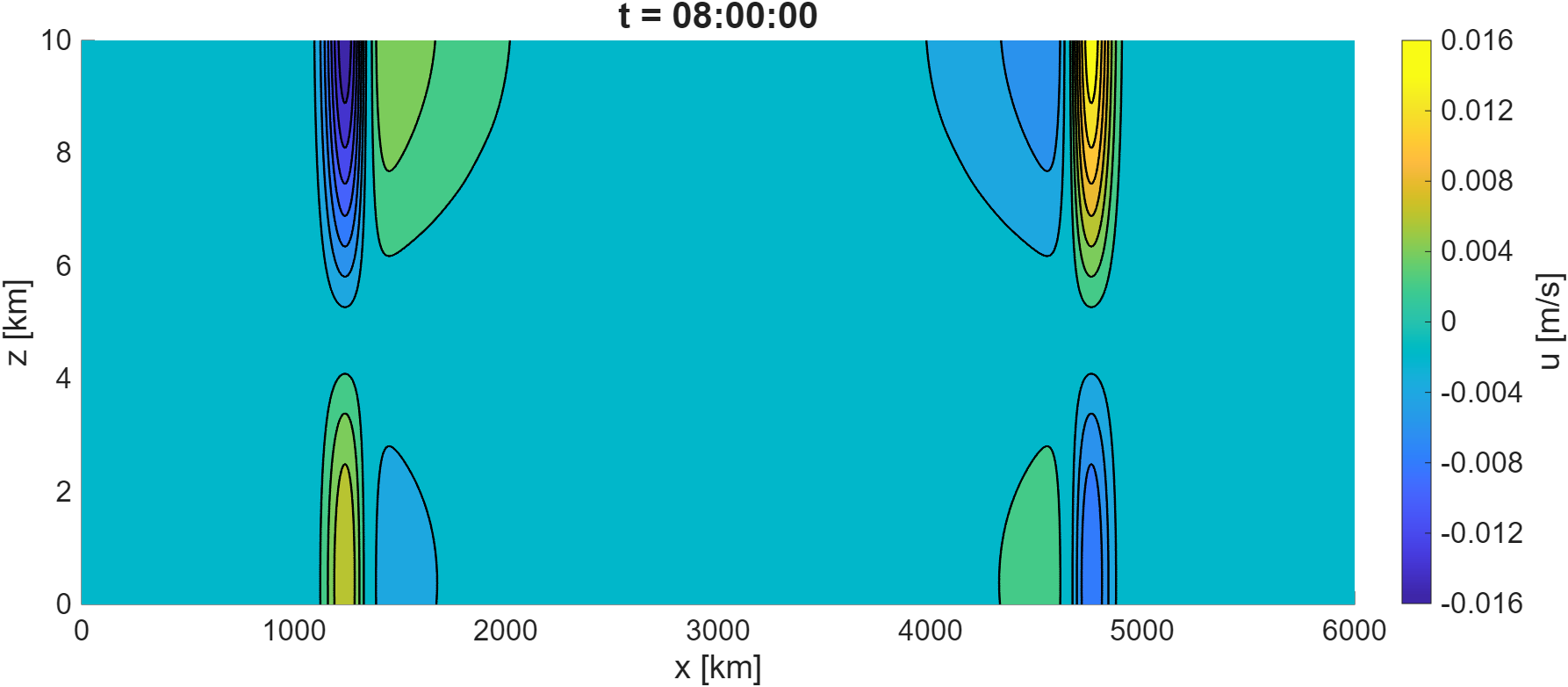}
    \caption{}
    \label{fig:Baldauf_u}
  \end{subfigure}

  \caption{Reference solution from \cite{baldauf_analytic_2013} (contours) and numerical solution computed by the IMEX-DG scheme (colour shading) for the inertia-gravity wave test case at final time $T_f=\SI{8}{\hour}$. Panel (a): temperature perturbation, contours in the range $\SI[parse-numbers=false]{[-0.006,0.006]}{\kelvin}$ with a $\SI{0.001}{\kelvin}$ interval. Panel (b): vertical velocity perturbation, contours in the range $\SI[parse-numbers=false]{[-0.0016,0.0016]}{\meter\per\second}$ with a $\SI{0.0002}{\meter\per\second}$ interval. Panel (c): horizontal velocity perturbation, contours in the range $\SI[parse-numbers=false]{[-0.016,0.016]}{\meter\per\second}$ with a $\SI{0.002}{\meter\per\second}$ interval.}
  \label{fig:Baldauf_qualitative_all}
\end{figure}

For the two-dimensional tests, a vertical-slice approximation is adopted to retain the effect of rotation on the vertical structure of the flow while reducing the problem to a two-dimensional physical domain. The main benchmark considered is the propagation of inertia-gravity waves in a rectangular channel with no background flow and non-zero rotation, following the standard test cases proposed in \cite{baldauf_analytic_2013,benacchio_semi-implicit_2019}.
The computational domain of the first case is a two-dimensional channel of size $\SI[parse-numbers=false]{(0,6000)\times(0,10)}{\kilo\meter\squared}$. Periodic boundary conditions are prescribed in the horizontal direction, while wall boundary conditions are imposed at the top and bottom boundaries.
The flow is initialized with a stationary, hydrostatic, horizontally homogeneous, isothermal background state: 
\begin{equation}
    \begin{aligned}
        T_0, \quad p_0(z)=p_s e^{-\delta z}, \quad \rho_0=\rho_s e^{-\delta z}, \quad \delta=\frac{g}{RT_0},\quad
        p_s=\frac{g \rho_s}{\delta}, \quad \vel_0=\bm{0},
         \label{eq:baldauf_background}
    \end{aligned}
\end{equation}
with $p_s=\SI[parse-numbers=false]{10^5}{\pascal}$ and $T_0=\SI[parse-numbers=false]{250}{\kelvin}$.
Following \cite{benacchio_semi-implicit_2019,baldauf_analytic_2013}, the thermodynamic perturbations are initialized through Bretherton-transformed variables
\begin{equation}
    \begin{aligned}
        T'(\mathbf{r},t=0)=e^{\tfrac{1}{2}\delta z} \cdot T_b(\mathbf{r},t=0), \qquad
        p'(\mathbf{r},t=0)=0, \qquad \rho'(\mathbf{r},t=0)=e^{-\tfrac{1}{2}\delta z} \cdot \rho_b(\mathbf{r},t=0),
       \label{eq:baldauf_bubble}
    \end{aligned}
\end{equation}
where the Bretherton-transformed temperature and density are given by
\begin{equation}
    \begin{aligned}
        T_b(\mathbf{r},t=0)=\Delta T \cdot e^{- \frac{(x-x_c)^2}{a^2}}\cdot \sin \pi \frac{z}{H}, \qquad
        \rho_b(\mathbf{r},t=0)=\rho_s \left( - \dfrac{T_b}{T_0} \right).
    \end{aligned}
\end{equation}
The perturbation bubble is centered at $x_c=\SI{3000}{\kilo\meter}$, with $a=\SI{100}{\kilo\meter}$ and $H=\SI{10}{\kilo\meter}$.
The velocity perturbation is initially set to zero. Rotation is included through a constant Coriolis parameter $f\approx \SI[parse-numbers=false]{1.03126 \cdot 10^{-4}}{\per\second}$, corresponding to a latitude of $\phi=\ang{45}$. 
The final time is set to $T_f=8$ h. To reproduce the effective resolution used in the DUNE simulations reported in \cite{baldauf_analytic_2013}, we consider a mesh composed by $(300 \times 20)$ elements and a polynomial degree $r=4$, yielding an effective resolution of $\SI{5}{\kilo\meter}$ in $x$ and $\SI{125}{\meter}$ in $z$. The time step is taken equal to $\SI{0.5}{\second}$, corresponding to an advective Courant number $C_{\mathrm{adv}}\approx2.28 \cdot 10^{-6}$ and an acoustic Courant number $C\approx 0.044$. As discussed in \cite{baldauf_analytic_2013}, a sufficiently small time step is necessary to properly resolve the acoustic–gravity wave dynamics and avoid spurious fast-wave activity contaminating the balanced solution.
Following \cite{baldauf_analytic_2013}, the analytical reference solution is obtained by linearising the fully compressible equations around an isothermal hydrostatic state, applying a Bretherton transformation, and evolving the resulting inertial gravity wave modes. This yields an exact solution of the linearised system used as a benchmark for fully nonlinear simulations.
Using a reference implementation to evaluate the analytical solution for this configuration, Figures ~\ref{fig:Baldauf_T},~\ref{fig:Baldauf_w}, and ~\ref{fig:Baldauf_u} compare the numerical solution against the analytical reference fields at the final time. The results are in line with those of \cite{baldauf_analytic_2013}, and confirm that the present solver is able to accurately capture the propagation and structure of the inertia–gravity waves emanating from the centre of the domain, where a geostrophic mode remains present until the end of the simulation.

\subsubsection*{Convergence test and comparison between approaches R1 and R2}
A convergence study is performed to assess the accuracy of the numerical scheme and to compare the two approaches (R1 and R2) introduced in Section \ref{sec:numerical_framework}. Both strategies are validated against the analytical reference solution at the final time $T_f = \SI{8}{\hour}$, following \cite{baldauf_analytic_2013}.

The error is evaluated for the perturbation variables $w$, $T'$, and $p'$ using the discrete $l^2$ and $l^\infty$ norms
\begin{equation}
\begin{aligned}
l^2(\psi) &= \left(\frac{1}{N_x N_z}\sum_{i=1}^{N_x}\sum_{k=1}^{N_z}|\psi(x_i,z_k)-\psi_{num}(x_i,z_k)|^2\right)^{1/2},\,\\
l^\infty(\psi) &= \max_{i,k} |\psi(x_i,z_k)-\psi_{num}(x_i,z_k)|,
\end{aligned}
\end{equation}
and the experimental order of convergence (EOC) is defined as
\begin{equation}
EOC_p = \frac{\log\left(l^p(\Delta x_2)/l^p(\Delta x_1)\right)}{\log\left(\Delta x_2/\Delta x_1\right)},
\end{equation}

The reference analytical solution is evaluated on a fixed grid with $N_x = 6000$ and $N_z = 400$, while the numerical solution is interpolated onto this grid for error computation.
The initial discretization consists of $150 \times 10$ elements, which corresponds to a resolution of $\SI{40}{\kilo\meter}$ in the $x$ direction and of $\SI{1}{\kilo\meter}$ in the $z$ direction, with $\Delta t = 1$ s and polynomial degree $r=4$. The standard  hyperbolic scaling of space and time intervals is adopted, in which the mesh resolution and the time step are both concurrently halved. The expected convergence rate is $\min\{2, r+1\}$.

\begin{table}[H]
\centering
\fontsize{8}{9.5}\selectfont
\caption{Computed errors and experimental orders of convergence with the IMEX-DG simulation of the inertia gravity wave test case: $w$ variable at $t=\SI{8}{\hour}$, R1 vs R2 strategies for the rotation and gravity terms.}
\label{tab:conv_w}
\begin{tabular}{r r r
r r r r
r r r r}
\toprule
& & \multicolumn{4}{c}{R1} & \multicolumn{4}{c}{R2} \\
\cmidrule(lr){4-7} \cmidrule(lr){8-11}
$\Delta x  $[km]& $\Delta z$[m] & $\Delta t \,[s]$
& $l^2$ & EOC$_2$ & $l^\infty$ & EOC$_\infty$
& $l^2$ & EOC$_2$ & $l^\infty$ & EOC$_\infty$ \\
\midrule
40&1000  & 1
& $3.83\cdot10^{-3}$ & --   & $2.35\cdot10^{-2}$ & --
& $3.83\cdot10^{-3}$ & --   & $2.35\cdot10^{-2}$ & -- \\

20&500  & 0.5
& $1.07\cdot10^{-4}$ & 5.16 & $5.17\cdot10^{-4}$ & 5.51
& $1.16\cdot10^{-4}$ & 5.10 & $5.36\cdot10^{-4}$ & 5.45 \\

10&250  & 0.25
& $2.28\cdot10^{-5}$ & 2.24 & $1.12\cdot10^{-4}$ & 2.21
& $2.38\cdot10^{-5}$ & 2.23 & $1.16\cdot10^{-4}$ & 2.20 \\

5&125 & 0.125
& $2.86\cdot10^{-6}$ & 3.00 & $1.62\cdot10^{-5}$ & 2.79
& $3.08\cdot10^{-6}$ & 2.95 & $1.43\cdot10^{-5}$ & 3.02 \\
\bottomrule
\end{tabular}
\end{table}

\begin{table}[H]
\centering
\fontsize{8}{9.5}\selectfont
\caption{As in Table \ref{tab:conv_w}, but for the $p'$ variable.}
\label{tab:conv_p}
\begin{tabular}{r r r
r r r r
r r r r}
\toprule
& & \multicolumn{4}{c}{R1} & \multicolumn{4}{c}{R2} \\
\cmidrule(lr){4-7} \cmidrule(lr){8-11}
$\Delta x  $[km]& $\Delta z$[m] & $\Delta t [s]$
& $l^2$ & EOC$_2$ & $l^\infty$ & EOC$_\infty$
& $l^2$ & EOC$_2$ & $l^\infty$ & EOC$_\infty$ \\
\midrule
40&1000    & 1
& $7.11\cdot10^{-2}$ & --   & $4.83\cdot10^{-1}$ & --
& $7.29\cdot10^{-2}$ & --   & $4.98\cdot10^{-1}$ & -- \\

20&500    & 0.5
& $2.20\cdot10^{-2}$ & 1.69 & $1.57\cdot10^{-1}$ & 1.62
& $2.27\cdot10^{-2}$ & 1.68 & $1.60\cdot10^{-1}$ & 1.64 \\

10&250    & 0.25
& $5.39\cdot10^{-3}$ & 2.03 & $4.13\cdot10^{-2}$ & 1.93
& $5.60\cdot10^{-3}$ & 2.02 & $4.25\cdot10^{-2}$ & 1.91 \\

5&125   & 0.125
& $1.87\cdot10^{-3}$ & 1.53 & $1.48\cdot10^{-2}$ & 1.48
& $1.88\cdot10^{-3}$ & 1.57 & $1.48\cdot10^{-2}$ & 1.52 \\
\bottomrule
\end{tabular}
\end{table}

Values close to the expected ones are obtained for the convergence rates both with the R1 and the R2 approach in the two variables (Tables~\ref{tab:conv_w} and \ref{tab:conv_p}). More in detail, for $w$ both approaches achieve near-optimal convergence rates consistent with the expected theoretical order. For $p'$, a reduced but consistent convergence trend is observed. For some other variables, the convergence rate decreases at the finest refinement levels, suggesting error saturation associated with the evaluation of the analytical reference solution rather than with the numerical discretisation itself. This interpretation is supported by additional self-convergence tests (not shown), where the expected asymptotic behaviour is recovered. 
Overall, these results indicate that the two rotational treatments are numerically equivalent for the considered test case, both in terms of convergence properties and solution accuracy. The R1 approach will be adopted in the subsequent numerical experiments.

\subsection{Planetary-scale inertia-gravity wave}

To further assess the robustness of the proposed formulation in a more demanding regime, the inertia-gravity wave benchmark is repeated on a larger computational domain. The horizontal extent of the domain is increased by one order of magnitude, yielding a domain of size $\SI[parse-numbers=false]{(0,60\,000) \times (0,10)}{\kilo\meter\squared}$. Consistently with the reference configuration considered in the previous test, the horizontal semi-amplitude of the initial perturbation is also increased by a factor of 10, corresponding to $a=\SI{1000}{\kilo\meter}$. All the remaining physical and numerical parameters are kept unchanged, including the background stratification, the Coriolis parameter, the number of elements of the mesh, the polynomial degree, and the boundary conditions. The simulation is run up to a longer final time $T_f= 4 \text{ days}=\SI{96}{\hour}$ using a time step $\Delta t = \SI{10}{\second}$. The corresponding acoustic and advective Courant numbers are approximately $C=0.09$ and $C_{\text{adv}}=2.81\times 10^{-6}$, respectively.

The larger horizontal scale leads to inertia-gravity waves with longer wavelengths and slower propagation speeds, while rotational effects become comparatively more significant in the overall flow evolution. This configuration therefore provides a more stringent test for the preservation of balanced dynamics and for the long-time robustness of the IMEX-DG formulation.

\begin{figure}[h!]
  \centering
  \begin{subfigure}{0.65\textwidth}
    \centering
    \includegraphics[width=\textwidth]{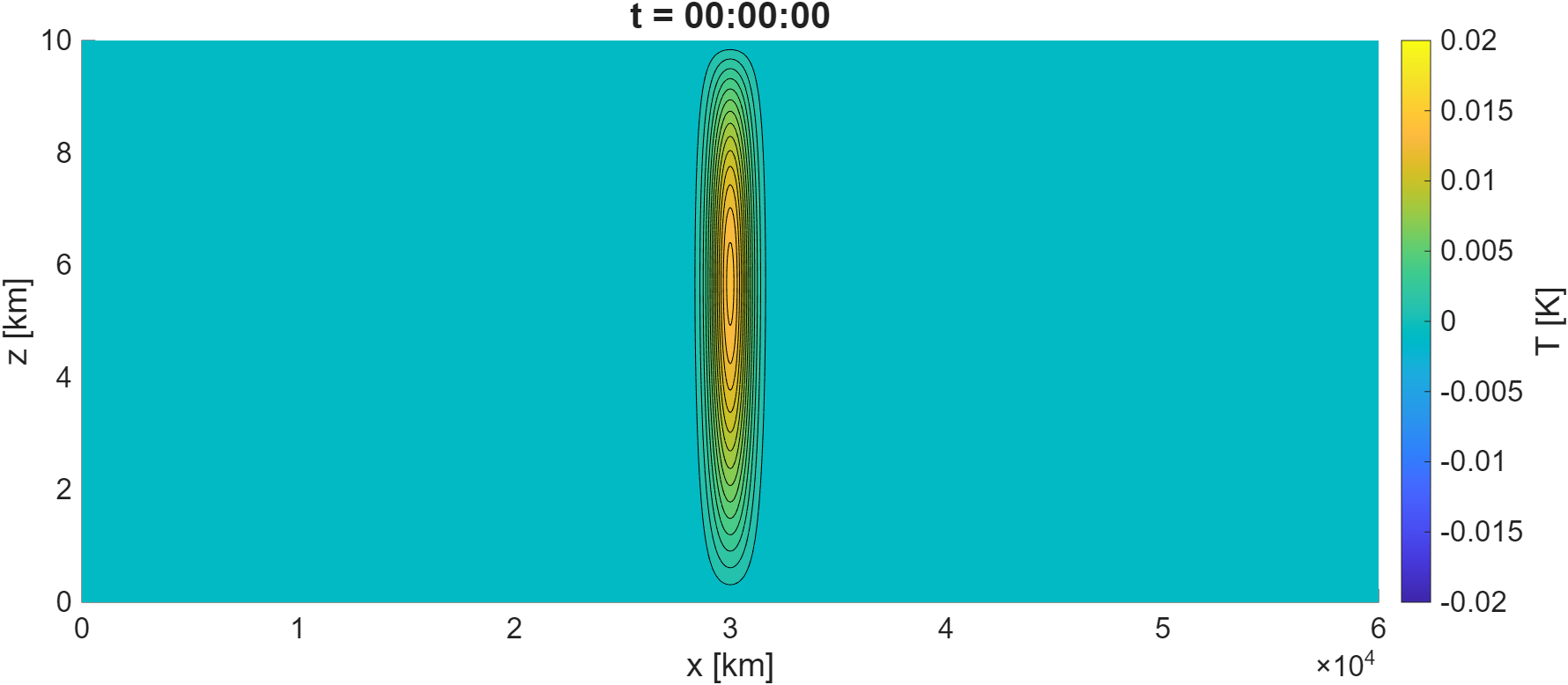}
    \caption{}
    \label{fig:IGWt0}
  \end{subfigure}
  \vspace{0.3cm}

  \begin{subfigure}{0.65\textwidth}
    \centering
    \includegraphics[width=\textwidth]{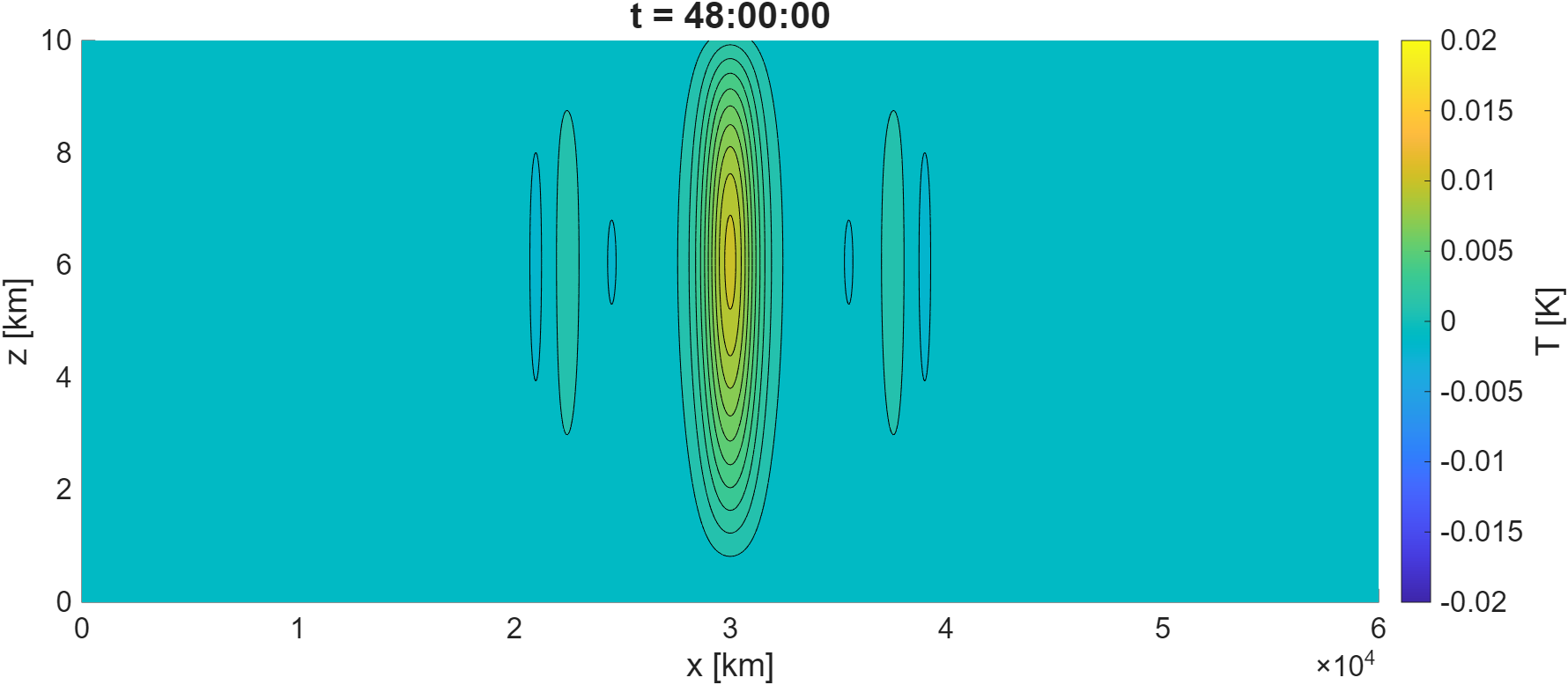}
    \caption{}
    \label{fig:IGWt2gg}
  \end{subfigure}

\vspace{0.3cm}

  \begin{subfigure}{0.65\textwidth}
    \centering
    \includegraphics[width=\textwidth]{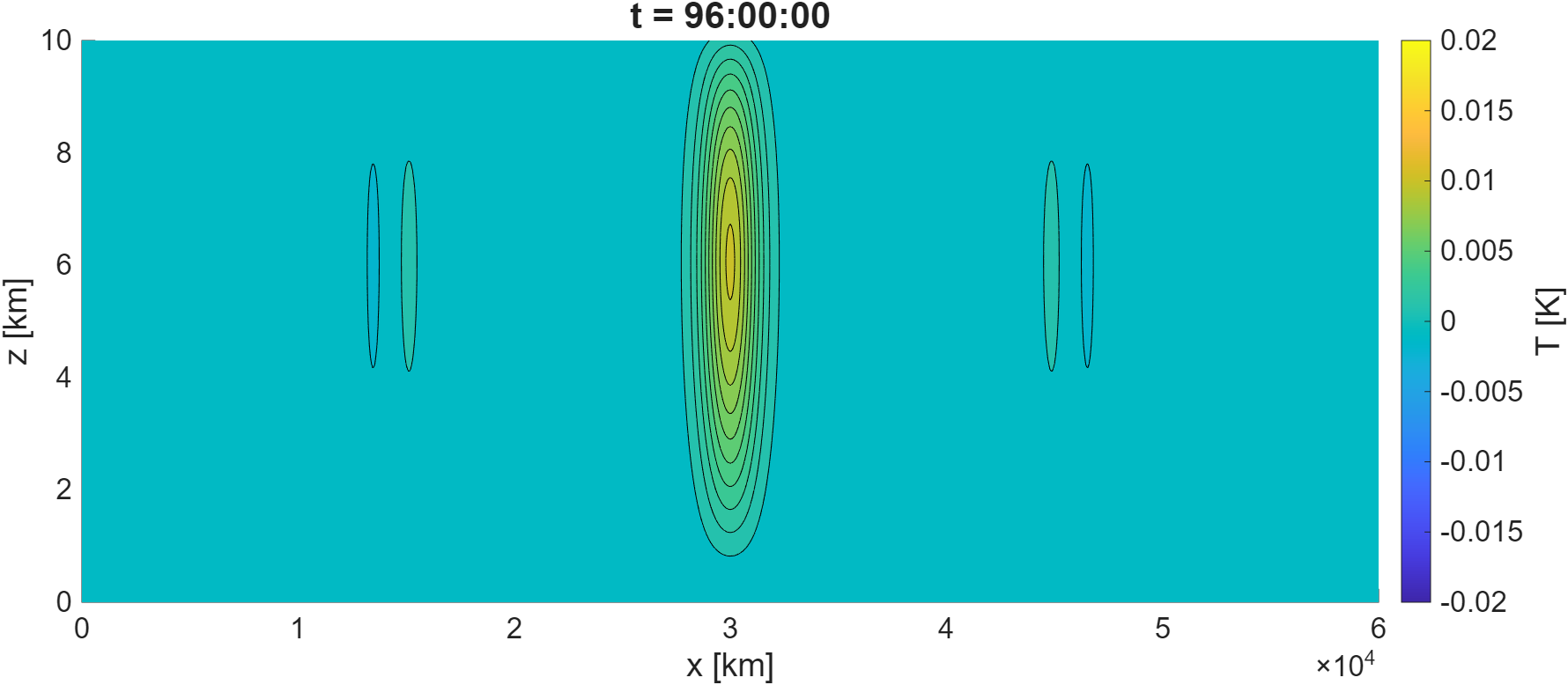}
    \caption{}
    \label{fig:IGW4gg}
  \end{subfigure}

  \caption{Temperature perturbation field for the planetary-scale inertia-gravity wave test at (a) $t=\SI{0}{\hour}$, (b) $t=\SI{48}{\hour}$, and (c) $t=\SI{96}{\hour}$. Contours from $\SI{-0.02}{\kelvin}$ to $\SI{0.02}{\kelvin}$ with an interval of $\SI{0.001}{\kelvin}$.}
  \label{fig:Baldauf_big_domain}
\end{figure}

Figure~\ref{fig:Baldauf_big_domain} shows the evolution of the temperature perturbation field at different times during the simulation. Initially, the perturbation is localized around the center of the domain and exhibits the expected vertically elongated structure associated with the imposed stratification. As the simulation evolves, the perturbation generates pairs of inertia-gravity waves propagating horizontally away from the central region. At $t=\SI{48}{\hour}$, the wave packets are clearly visible on both sides of the initial disturbance while the central balanced structure remains well preserved. At the final time $t=\SI{96}{\hour}$, the waves have propagated further across the domain without exhibiting noticeable numerical distortion or spurious oscillations.

Overall, the numerical solution remains stable over the entire simulation and preserves the expected symmetry and large-scale wave structure. These results confirm that the rotating IMEX-DG formulation remains robust and accurate also in large-domain configurations characterized by long integration times and stronger rotational influence.

To further assess the ability of the scheme to preserve balanced dynamics over long integration times, we also provide a quantitative measure of the departure from geostrophic equilibrium during the simulation. Since the present formulation is not exactly well-balanced at the discrete level, the pressure-gradient and Coriolis contributions do not cancel to machine precision, and small imbalance errors may progressively develop during the time integration.
The discrete residual of the horizontal momentum balance is computed, following \cite{tumolo_semi-implicit_2013}, as
\begin{equation}
\label{eq:egeo}
\mathcal{E}_{\mathrm{geo}}
= \frac{\int_0^{z_{\text{max}}}\int_0^{x_{\text{max}}}\left|\frac{\partial p}{\partial x}\frac{1}{\rho}-f v\right| \mathrm{d}x \mathrm{d}z}{z_{\text{max}}x_{\text{max}}},
\end{equation}
where $x_{\text{max}}$ and $z_{\text{max}}$ denote the the maximum horizontal and vertical coordinates of the domain, respectively. The integrals are evaluated numerically using the trapezoidal rule.

\begin{figure}[H]
  \centering
  \includegraphics[width=0.5\linewidth]{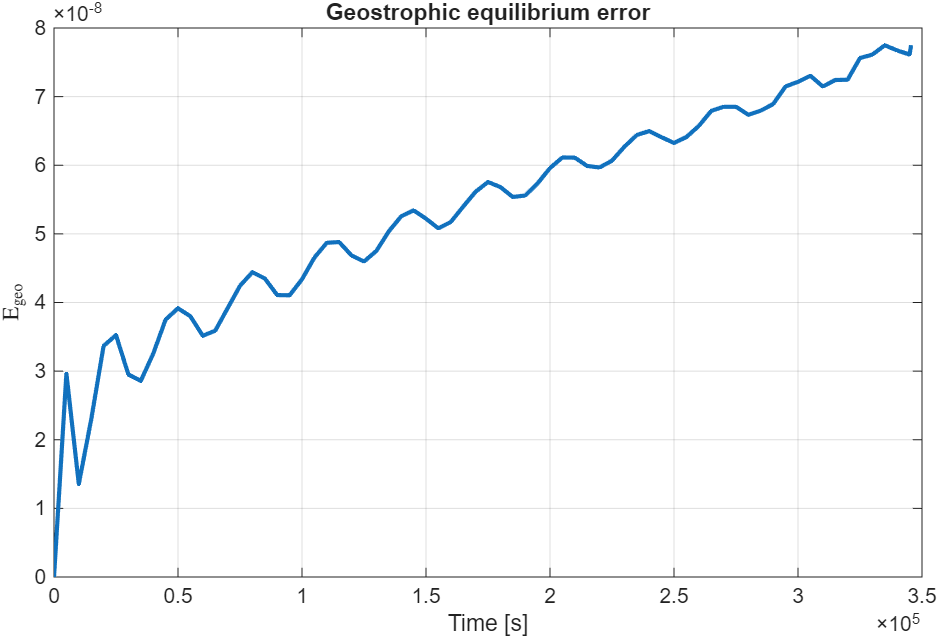}
  \caption{Time evolution of the discrete geostrophic equilibrium error $\mathcal{E}_{\mathrm{geo}}$ during the planetary-scale inertia-gravity wave simulation.}
  \label{fig:geo_error_large_domain}
\end{figure}

Although a slow growth of $\mathcal{E}_{\mathrm{geo}}$ is observed over the four-day simulation time, the magnitude of the error remains lower than $10^{-7}$ (Figure~\ref{fig:geo_error_large_domain}), which is compatible with the truncation error of the numerical discretisation. Moreover, the error growth remains bounded over the full integration interval and does not trigger spurious oscillations or visible distortions of the large-scale wave structures. This behaviour indicates that, despite the lack of a geostrophically and hydrostatically well-balanced discretisation, the proposed high-order IMEX-DG formulation preserves geostrophic equilibrium with satisfactory accuracy and robustness also over long simulation times and in large-scale rotating configurations \cite{baldauf_comparison_2026}. Recent results with DG discretisations of the compressible Euler equations show that the choice of the prognostic variable and the resulting discrete formulation can affect the preservation of hydrostatic equilibria \cite{baldauf_comparison_2026}. In this context, the present results suggest that the high-order discretisation keeps the geostrophic imbalance small in the considered regime, despite the absence of an exactly well-balanced formulation. Future work will therefore investigate more systematically the regimes in which such a formulation can provide sufficiently small equilibrium imbalances for longer times, through a comparison with a well-balanced formulation. 

\subsection{3D flow over a bell-shaped hill}
\label{sec:3D}
In this section, we consider fully three-dimensional atmospheric dynamics simulations with the IMEX-DG scheme. These experiments enable the assessment of the interaction between stratification, orographic forcing, and rotation in more realistic and representative dynamical configurations. 
Two test cases of flow over orography are considered: an existing medium-steep benchmark \cite{melvin_mixed_2019, orlando_robust_2024, orlando_efficient_2025, orlando_improving_2025}, used  to analyse the parallel performance of the solver, and a new, larger-scale atmospheric configuration designed to investigate the behaviour of the method in a regime where the effect of rotation is more prominent. Both configurations are simulated with and without rotation in order to assess the impact of the Coriolis term and verify the robustness of the solver in fully three-dimensional flows.
For both tests, the computational domain is a hexahedral box where the lower boundary is deformed by a smooth bell-shaped mountain centered at $(x_c,y_c)$ km in the horizontal plane, with height 
\begin{equation}
    h(x,y)=h_0\left(1+\frac{(x-x_c)^2+(y-y_c)^2}{a_c^2}\right)^{-1},
    \label{eq:bell_mountain}
\end{equation}
where $h_0$ denotes the maximum mountain height and $a_c$ controls the horizontal half-width \cite{orlando_robust_2024}.
The mesh is generated on the physical domain and mapped to the terrain-following coordinates through the Gal-Chen transformation \cite{gal-chen_use_1975}. Periodic boundary conditions are imposed in the horizontal directions, while solid-wall conditions are applied at the bottom boundary. At the upper boundary, non-reflecting conditions are employed in order to reduce artificial wave reflections. As in \cite{orlando_impact_2024,orlando_robust_2024}, to further suppress spurious reflections of outgoing gravity waves, Rayleigh damping (sponge) layers are introduced near the upper boundary and in lateral buffer regions. The damping coefficient is set to zero below a prescribed height $z_B$ and increases smoothly toward the top boundary $z_T$ according to
\begin{equation}
\begin{aligned}
    &\lambda(z)=\bar{\lambda}\sin^2\left[\frac{\pi}{2}\left(\frac{z-z_B}{z_T-z_B}\right)\right]
    \qquad &&\text{for } z\ge z_B, \\
    &\lambda(z)=0
    \qquad &&\text{for } z<z_B .
\end{aligned}
\label{eq:damping}
\end{equation}
An analogous formulation is employed in the lateral damping regions.

\subsubsection*{Small domain test}
The first domain considered is a $\Omega=\SI[parse-numbers=false]{(0,60) \times (0,40) \times (0,16)}{\kilo\meter\cubed}$ box. The center of the mountain is at $(x_c,y_c)=\SI[parse-numbers=false]{(30,20)}{\kilo\meter}$ and $h_0=\SI{400}{\meter}$, $a_c=\SI{1000}{\meter}$. Damping is applied in the first and last $\SI{20}{\kilo\meter}$ in the $x$ direction, in the first and last $\SI{10}{\kilo\meter}$ in the $y$ direction, and in the uppermost $\SI{6}{\kilo\meter}$ in the $z$ direction, with $\bar{\lambda} \Delta t=1.2$. The flow is initialized from a stationary hydrostatic background state. The expression for the background potential temperature is given by $\bar \theta=\theta_{\mathrm{ref}}\exp(\frac{N^2}{g}z)$. The buoyancy frequency is set to $N=\SI{0.01}{\per\second}$, corresponding to a nonhydrostatic regime, and $\theta_{\mathrm{ref}}=\SI{293.15}{\kelvin}$.
When rotation is active, the background is completed with a geostrophically balanced horizontal pressure gradient. The background pressure is decomposed into a height-dependent hydrostatic contribution and a geostrophic correction accounting for rotation:
\begin{equation}
\label{eq:init_cond}
\begin{aligned}
     &p_0=p_{0,h}(z)+p_{0,g}(x,y), \\
   & p_{0,h}= p_{\text{ref}}(\pi)^{\frac{\gamma}{\gamma-1}} ,\quad  \pi =1+ \frac{g^2}{c_p T_{\text{ref}}N^2} \left[\exp{\left(-\frac{N^2}{g} z\right)}-1\right],\\
   & p_{0,g}(x,y)= (x-x_0)\,\rho_0 f\, U\sin\alpha-(y-y_0)\,\rho_0 f\, U\cos\alpha,
\end{aligned}
\end{equation}
where $\pi$ is the Exner pressure. In the non-rotating case ($f=0$), the geostrophic contribution vanishes and $p_0$ reduces to the hydrostatic profile. A uniform horizontal wind is imposed,  $\mathbf{u}_0=U_0(\cos{\alpha},\sin{\alpha},0)^{\top},$ with $\alpha=0$ (flow aligned with the $x$ direction) and $U_0=\SI{10}{\meter\per\second}$. 
We use polynomial degree $r=4$ on a uniform mesh of $30 \times 20 \times 8 = 4800$ elements, yielding an effective resolution of $\SI{500}{\meter}$ in all directions, and $6 \times 10^5$ degrees of freedom for each scalar variable. The time step is $\SI{8}{\second}$ and the final time is $T_f=8$ h.  The corresponding acoustic and advective Courant numbers are approximately $C=5.5$ and $C_{\text{adv}}=0.21$, respectively. We consider two simulations: a non-rotating run ($f=0$) and a rotating run with $f=\SI[parse-numbers=false]{10^{-4}}{\per\second}$.

\begin{figure}[H] 
\centering

\begin{subfigure}[t]{0.48\textwidth}
  \centering
  \includegraphics[width=\linewidth]{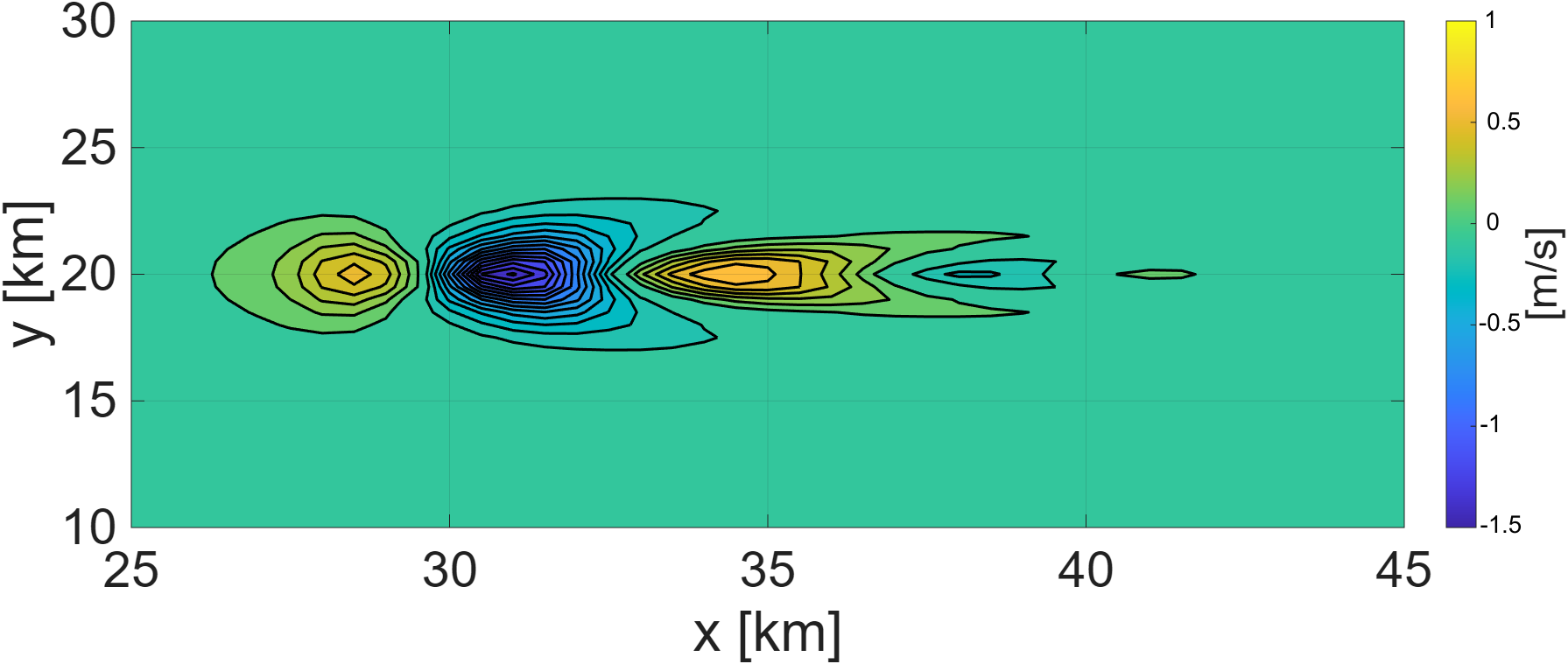}
  \caption{}
\end{subfigure}\hfill
\begin{subfigure}[t]{0.48\textwidth}
  \centering
  \includegraphics[width=\linewidth]{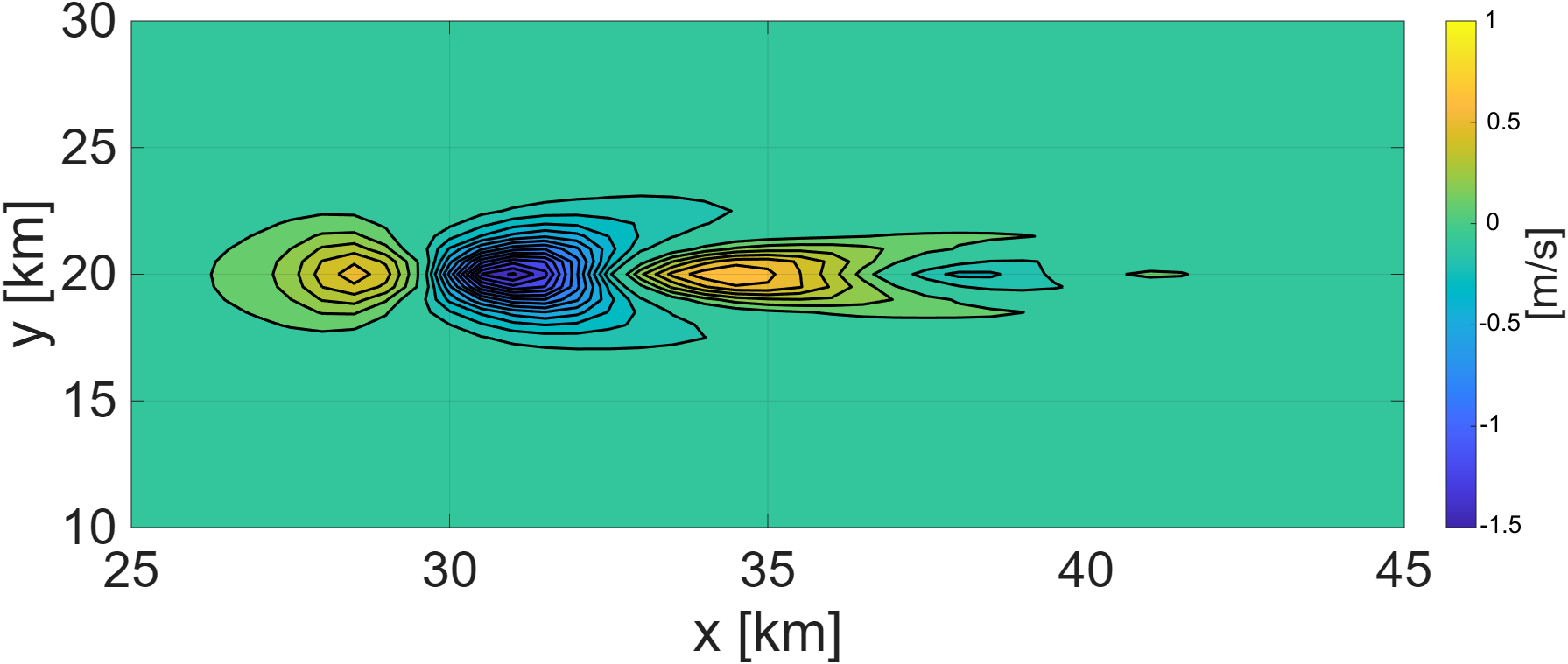}
  \caption{}
\end{subfigure}

\vspace{0.5em}

\begin{subfigure}[t]{0.48\textwidth}
  \centering
  \includegraphics[width=\linewidth]{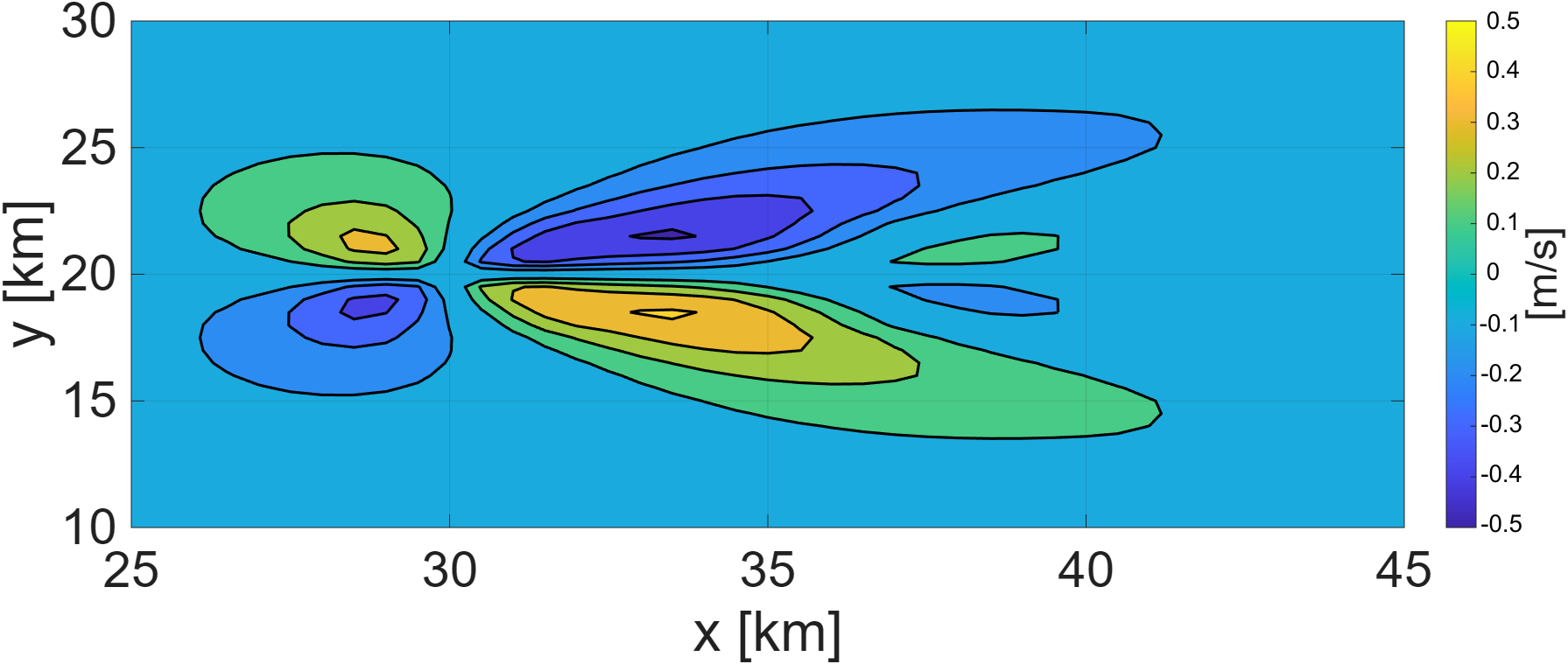}
  \caption{}
\end{subfigure}\hfill
\begin{subfigure}[t]{0.48\textwidth}
  \centering
  \includegraphics[width=\linewidth]{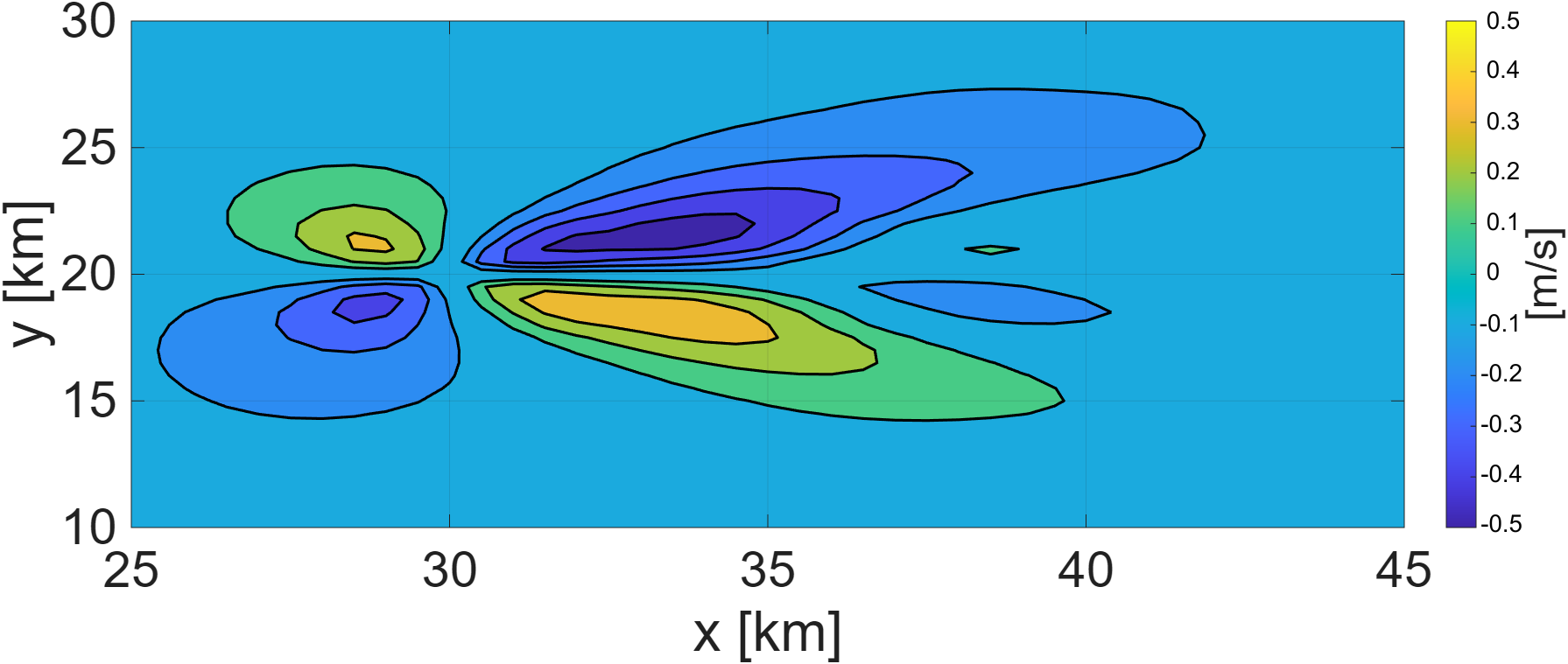}
  \caption{}
\end{subfigure}

\vspace{0.5em}

\begin{subfigure}[t]{0.48\textwidth}
  \centering
  \includegraphics[width=\linewidth]{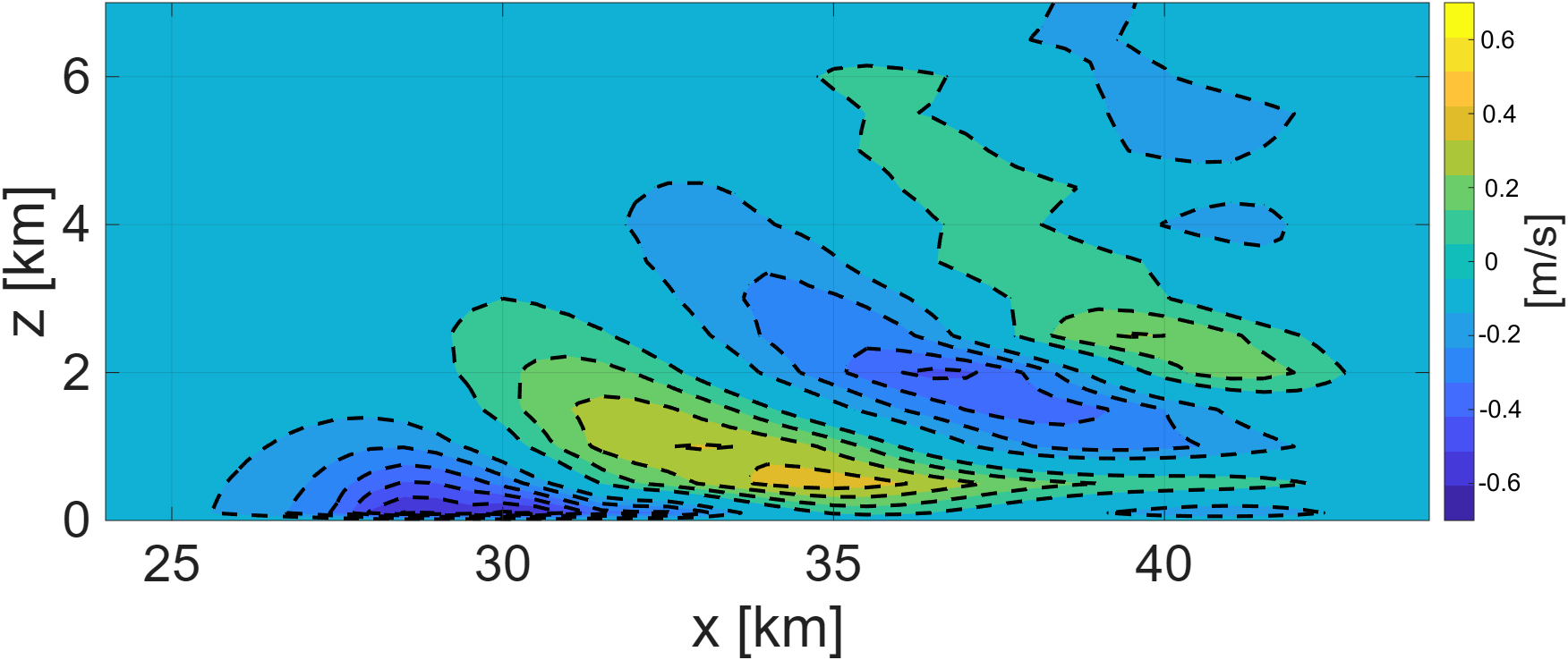}
  \caption{}
\end{subfigure}\hfill
\begin{subfigure}[t]{0.48\textwidth}
  \centering
  \includegraphics[width=\linewidth]{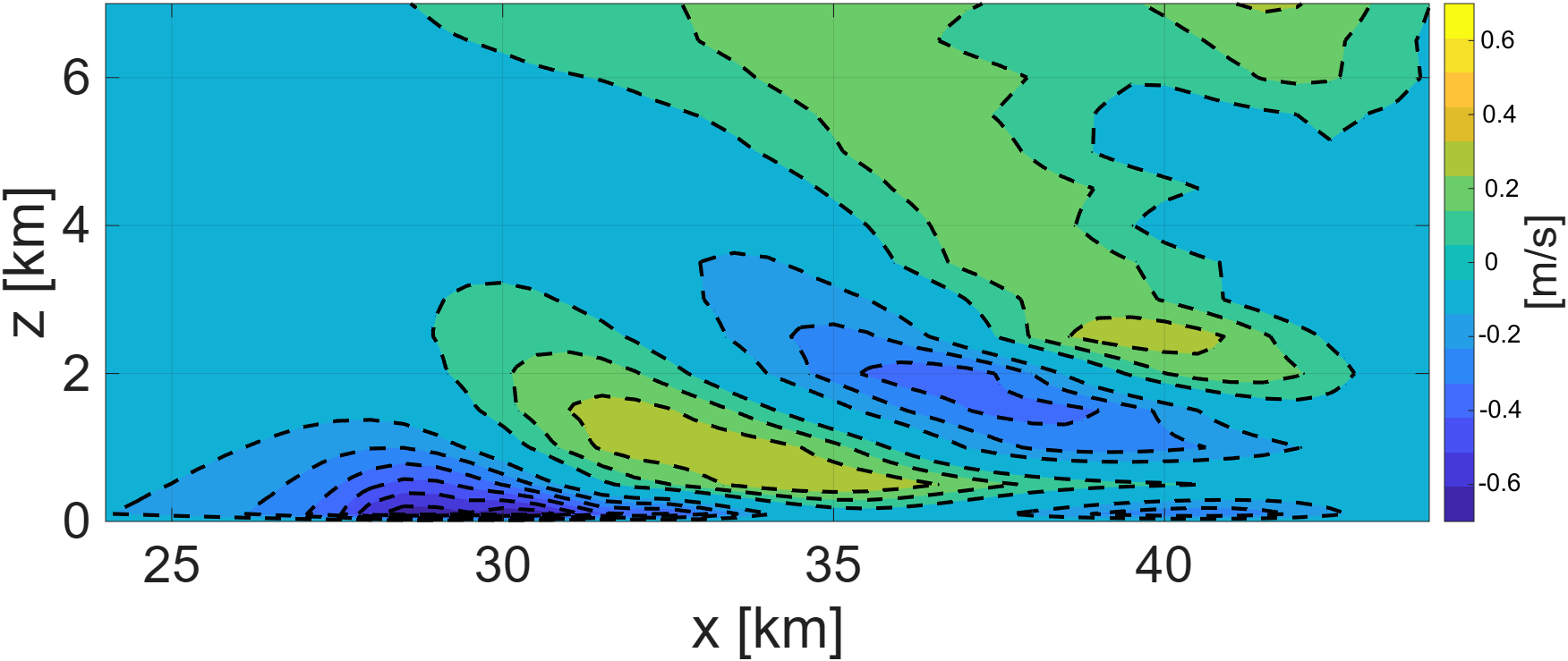}
  \caption{}
\end{subfigure}
  \begin{subfigure}[t]{0.48\textwidth}
  \centering
  \includegraphics[width=\linewidth]{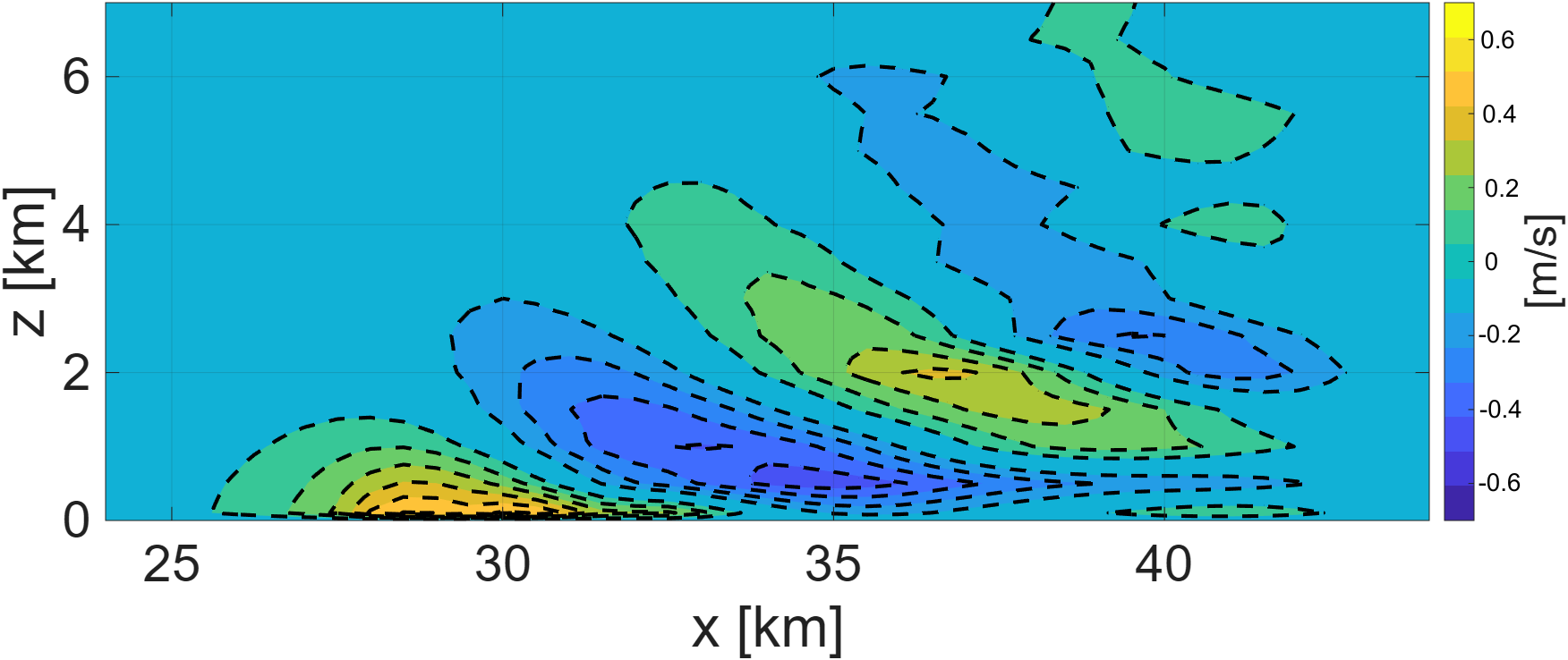}
  \caption{}
\end{subfigure}\hfill
\begin{subfigure}[t]{0.48\textwidth}
  \centering
  \includegraphics[width=\linewidth]{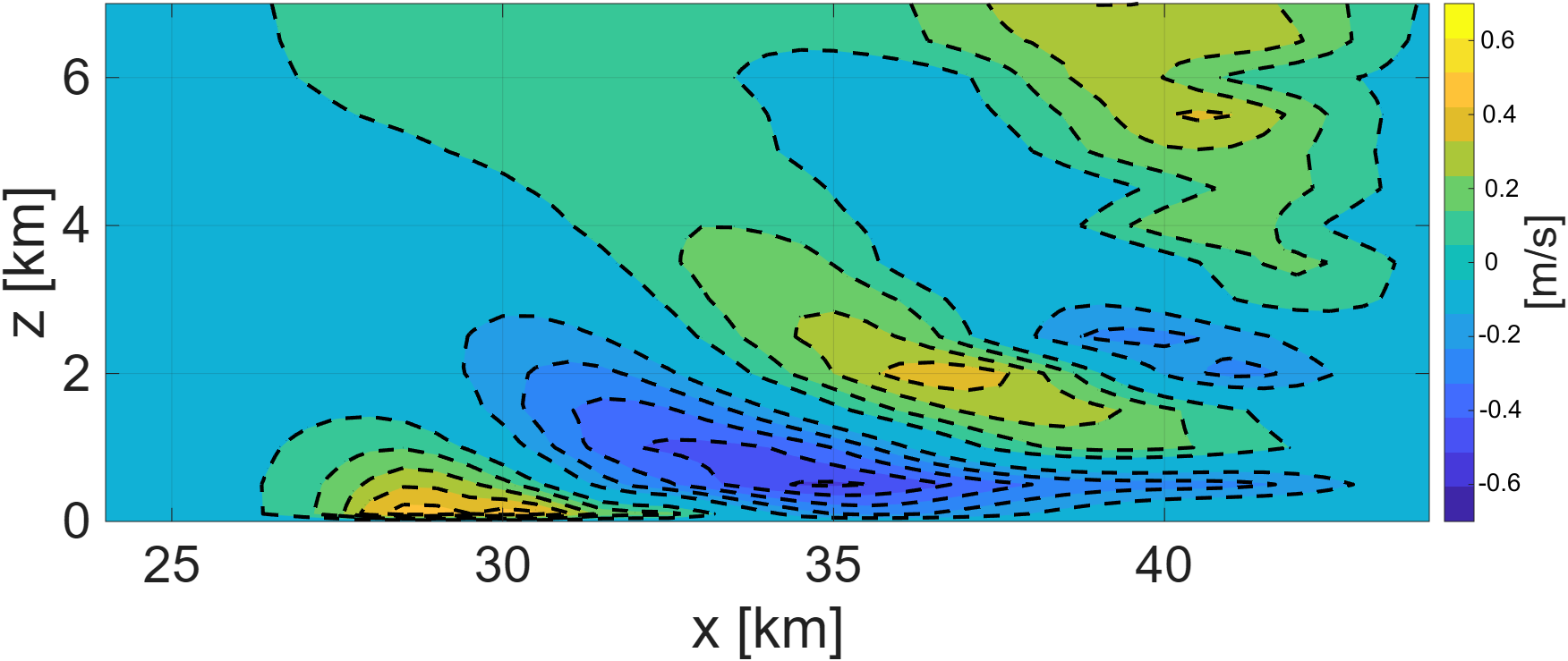}
  \caption{}
\end{subfigure}

\caption{Computed velocity at final time $T_f=\SI{8}{\hour}$ for the small domain test, non-rotating (left column) and rotating (right column) flow over a 3D bell-shaped hill. Panels (a) and (b) show the $u$ component on a $x-y$ slice at $z=\SI{800}{\meter}$, contours in the range $\SI[parse-numbers=false]{[-1.5,1.3]}{\meter\per\second}$ with a $\SI{0.1}{\meter\per\second}$ interval. Panels (c) and (d) show the $v$ component on a $x-y$ slice at $z=\SI{800}{\meter}$, contours in the range $\SI[parse-numbers=false]{[-0.5,0.5]}{\meter\per\second}$ with a $\SI{0.1}{\meter\per\second}$ interval. The bottom four panels show the $v$ component on a $x-z$ slice at $y=\SI{18}{\kilo\meter}$ (panels (e) and (f)), and at $y=\SI{22}{\kilo\meter}$ (panels (g) and (h)), contours in the range $\SI[parse-numbers=false]{[-0.7,0.7]}{\meter\per\second}$ with a $\SI{0.1}{\meter\per\second}$ interval. }
\label{fig:small_domain_comparison}
\end{figure} 

In both configurations, the orographic forcing generates a wake and a downstream gravity-wave pattern with comparable peak amplitudes, indicating that the inclusion of the Coriolis terms does not introduce unphysical effects (Figure \ref{fig:small_domain_comparison}). However, in the rotating case, the downstream structure loses the axial symmetry of the non-rotating case, consistently with the effect of rotation. This is evident in the horizontal slices of the $v$ component (panel (d) in Figure \ref{fig:small_domain_comparison}), and even more clearly in the vertical sections (third and fourth row in Figure \ref{fig:small_domain_comparison}): at the same downstream distance from the mountain, the non-rotating case exhibits contours that are mirror images with opposite sign across the wake, whereas the rotating case clearly shows a difference in the contour shapes, revealing the induced asymmetry. Since the computational domain is relatively limited in horizontal extent, the rotational effects remain moderate. Overall, this benchmark provides a first validation of the implementation of the rotating terms in the IMEX-DG scheme in a fully three-dimensional setting.

To obtain cleaner diagnostics while keeping the computational cost manageable, two additional simulations are performed using a non-conforming mesh. The non-conforming discretization allows for local refinement near the orography, where the strongest dynamics and wave generation occur, while maintaining a coarser resolution elsewhere \cite{orlando_robust_2024}.

\begin{figure}[H]
    \centering
    \includegraphics[width=0.6\linewidth]{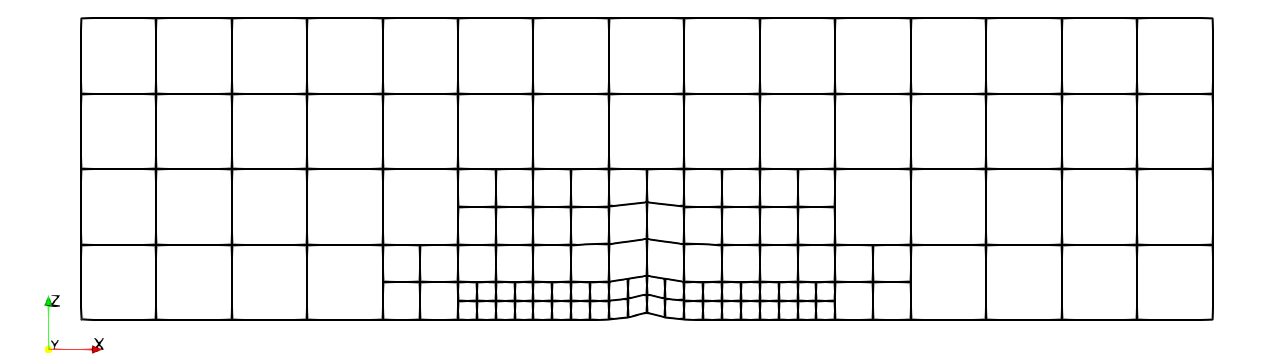}
    \caption{Section of the non-conforming mesh at $y=20$ km for the small domain test.}
    \label{fig:mesh_non_conf}
\end{figure}

The base mesh has an effective horizontal resolution of $1$ km and is locally refined by two additional levels around the mountain, reaching a finest resolution of approximately $250$ m. Refinement is applied in the region $x\in[20,40]$ km and $y\in[10,30]$ km; vertically, the first refinement level is used for $z\le7$ km and the second for $z\le2$ km (Figure \ref{fig:mesh_non_conf}), following the same approach as in \cite{orlando_impact_2024, orlando_robust_2024}. 

With this non-conforming mesh, the total number of active cells is 1594 (about one third of the conforming mesh), and the number of DOF per scalar variable is $199\,250$, i.e. also about one third of the uniform-grid configuration.

\begin{figure}[H] 
\centering
\begin{subfigure}[t]{0.48\textwidth}
  \centering
  \includegraphics[width=\linewidth]{figures/rot_vertical_18.png}
  \caption{}
\end{subfigure}\hfill
\begin{subfigure}[t]{0.48\textwidth}
  \centering
  \includegraphics[width=\linewidth]{figures/rot_vertical_22.png}
  \caption{}
\end{subfigure}

\vspace{0.6em}

\begin{subfigure}[t]{0.48\textwidth}
  \centering
  \includegraphics[width=\linewidth]{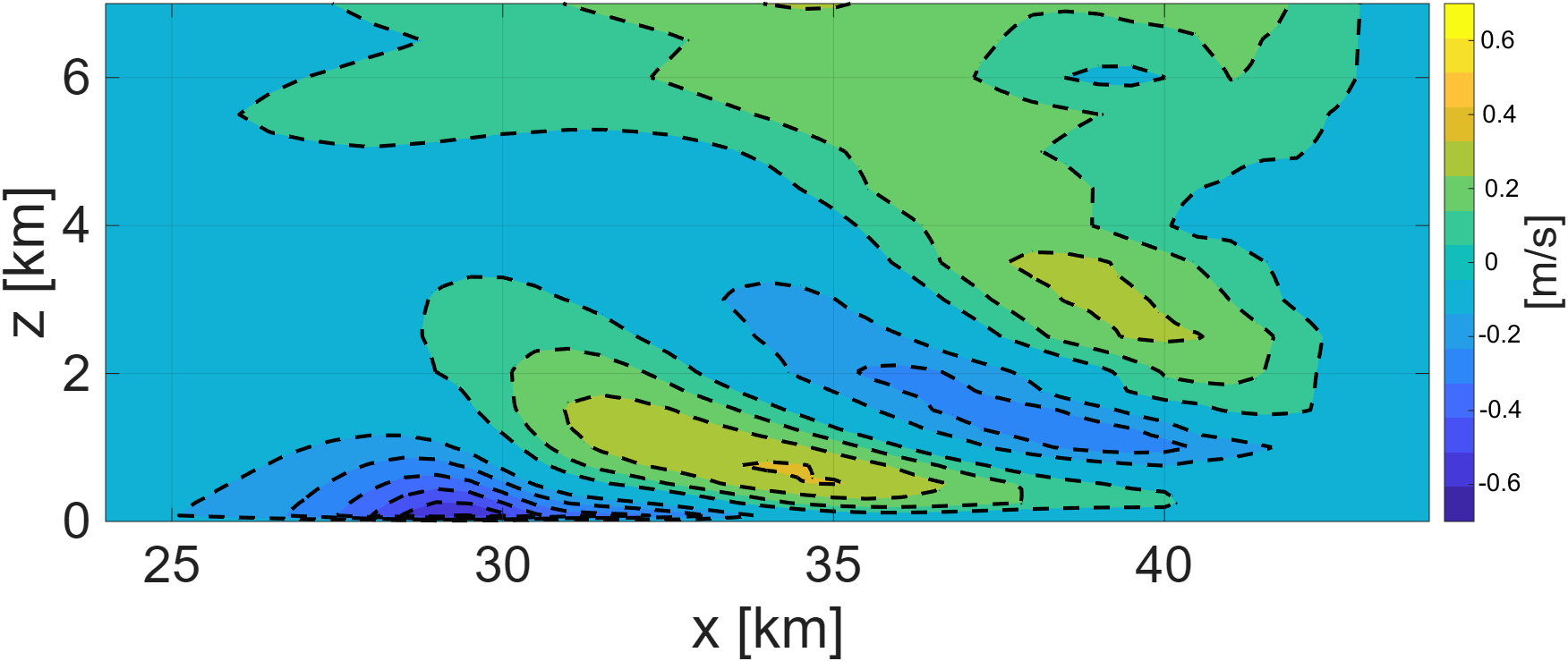}
  \caption{}
\end{subfigure}\hfill
\begin{subfigure}[t]{0.48\textwidth}
  \centering
  \includegraphics[width=\linewidth]{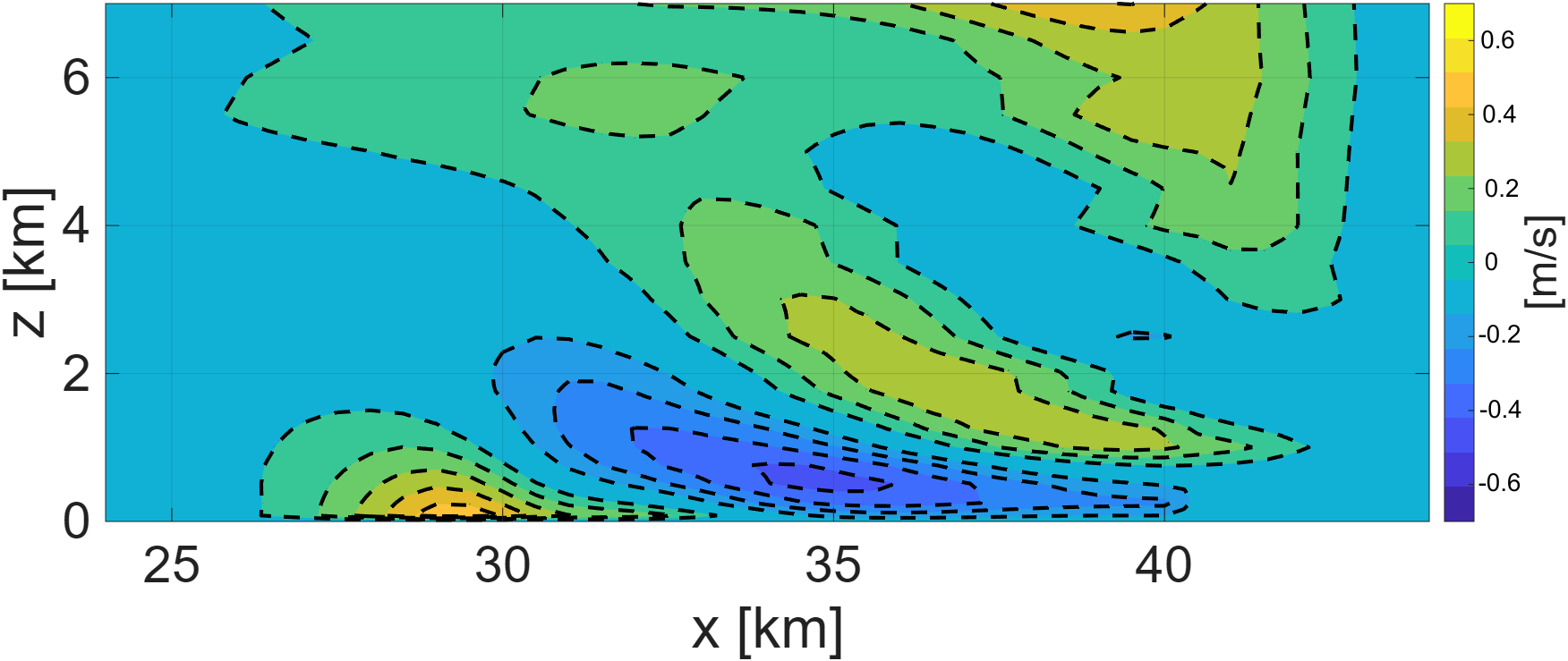}
  \caption{}
\end{subfigure}
\caption{Small domain test: comparison of the computed $v$ velocity component at final time $T_f=\SI{8}{\hour}$ using a uniform mesh (first row) and non-conforming mesh (second row). Vertical $x-z$ slices at $y=\SI{18}{\kilo\meter}$ (left column) and at $y=\SI{22}{\kilo\meter}$ (right column), contours in the range $\SI[parse-numbers=false]{[-0.7,0.7]}{\meter\per\second}$ with a $\SI{0.1}{\meter\per\second}$ interval.}
\label{fig:small_domain_nonconf}
\end{figure}

The simulation using the non-conforming mesh reproduces the same large-scale wave dynamics observed with the uniform grid, while providing smoother and better-resolved structures near the mountain and in the near-wake region, where the strongest gradients occur (Figure \ref{fig:small_domain_nonconf}). The comparison between rotating and non-rotating cases leads to the same qualitative conclusions as before, with rotation producing only a weak asymmetry in the downstream wave pattern. Overall, the non-conforming discretization improves the local accuracy around the orography at a significantly reduced computational cost, while also enabling targeted vertical refinement near the terrain.

\begin{figure}[H]
  \centering
 \includegraphics[width=0.6\textwidth]{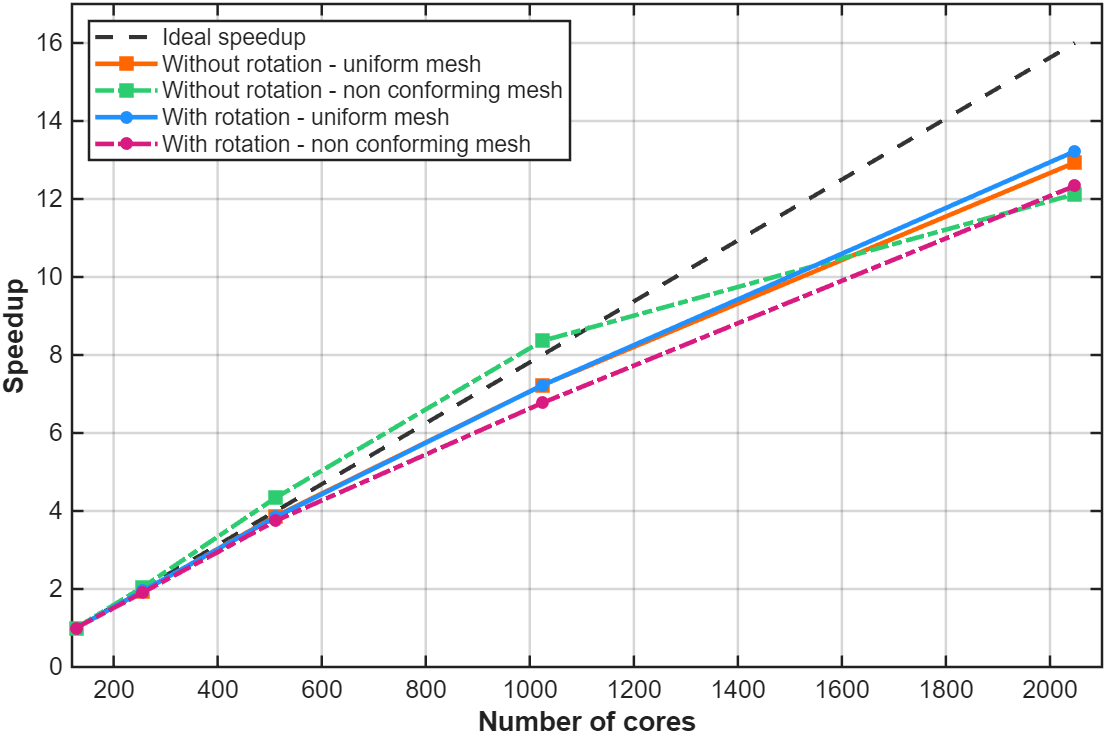}
  \caption{Small domain test: strong-scaling results for the rotating (circle markers) and the non-rotating (square markers) test cases using both uniform (solid lines) and non-conforming (dashed-dotted lines) meshes. The non-rotating results correspond to those in \cite{orlando_improving_2025}} 
  \label{fig:strong_scaling}
\end{figure}

Since this benchmark coincides with the configuration considered in \cite{orlando_improving_2025,orlando_efficient_2025} for the non-rotating case, we repeat the strong-scaling analysis in the presence of rotation in order to assess the impact of the rotational term on the parallel performance of the solver and to directly compare the results with the previously published results.
We consider both a uniform mesh composed of $120\times80\times32 = 307\,200$ elements with polynomial degree $r=4$, and a non-conforming mesh consisting of $204\,816$ elements.
The strong-scaling analysis evaluates the wall-clock time required to perform the simulations while keeping the computational workload, namely the spatial resolution, fixed, and progressively increasing the available computational resources. In particular, the simulations are executed using 1, 2, 4, 8, and 16 full MeluXina CPU nodes, each providing 128 cores. Figure \ref{fig:strong_scaling} reports the comparison between the previously published strong-scaling results obtained in the absence of rotation and the corresponding results that include the rotational term.

Overall, the introduction of rotation does not significantly deteriorate the parallel efficiency of the solver. The speedup curves obtained with and without rotation exhibit a very similar trend for both the uniform and non-conforming meshes, confirming that the additional rotational contribution introduces only a limited overhead in the parallel implementation. In particular, the solver maintains a close-to-ideal scaling behaviour up to 1024 cores, while a moderate degradation is observed at 2048 cores, where communication costs and load-balancing effects become more relevant. Moreover, the non-conforming configuration shows a scaling behaviour comparable to the uniform mesh case, demonstrating the robustness of the adaptive mesh refinement (AMR) framework also in the presence of rotational effects.

\subsubsection*{Large domain test}

In order to better assess the impact of rotation on the flow dynamics and to verify the correct preservation of geostrophic balance by the IMEX-DG solver, we designed a new test of 3D flow over orography on a larger computational domain. The enlarged domain allows rotational effects to develop over larger spatial scales, reducing the influence of the lateral damping regions and making the differences between the rotating and non-rotating configurations more evident.
The computational domain is a $\Omega=\SI[parse-numbers=false]{(0,240)\times(0,240)\times(0,30)}{\kilo\meter\cubed}$ box. The mountain is centered at $(x_c,y_c)=\SI[parse-numbers=false]{(120,120)}{\kilo\meter}$ and is defined by \eqref{eq:bell_mountain} with $h_0=\SI{100}{\meter}$ and $a_c=\SI{10}{\kilo\meter}$. Damping layers are applied in the first and last $\SI{60}{\kilo\meter}$ of both horizontal directions and in the uppermost $\SI{14}{\kilo\meter}$ of the domain, using expressions \eqref{eq:damping} with $\bar{\lambda}\Delta t = 1.2$.\\
The flow is initialized from a stationary hydrostatic background state. In this configuration, the Exner pressure profile is prescribed as
\begin{equation}
\pi=\left(\frac{p_0}{\bar{p}}\right)^{\frac{\gamma-1}{\gamma}}
=\exp\left(-\frac{g}{c_p\bar{T}}z\right),
\end{equation}
corresponding to an isothermal hydrostatic atmosphere. As in the previous test case, a uniform horizontal background wind is imposed, $\mathbf{u}_0=U_0(\cos\alpha,\sin\alpha,0)^{\top}$, with $\alpha=0$ and $U_0=\SI{10}{\meter\per\second}$. The reference potential temperature is set to $\theta_{\mathrm{ref}}=\SI{250}{\kelvin}$. The buoyancy frequency is equal to $N=\SI{0.02}{\per\second}$. The resulting non-dimensional parameter $\frac{Na_c}{U_0}=20\gg1$, so, according to the classification in \cite{pinty_simple_1995}, this configuration falls in the hydrostatic regime.

The simulations are performed using the polynomial degree $r=4$ on a uniform mesh composed of $50 \times 50 \times 50$ elements. The time step is $\Delta t=\SI{4}{\second}$ and the final time is set to $T_f=\SI{8}{\hour}$. The corresponding acoustic and advective Courant numbers are approximately $C=1.29$ and $C_{\text{adv}}=0.04$, respectively. 

\begin{figure}[H] 
\centering

\begin{subfigure}[t]{0.48\textwidth}
  \centering
  \includegraphics[width=\linewidth]{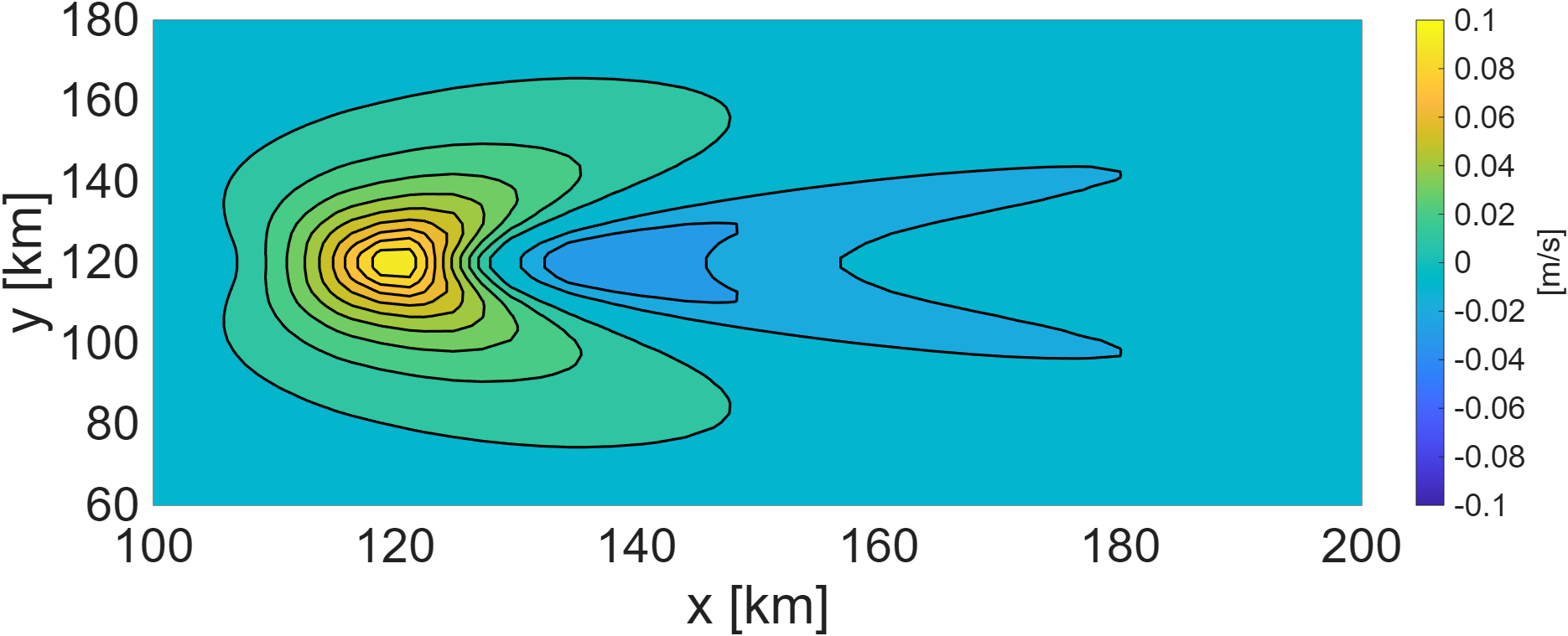}
  \caption{}
\end{subfigure}\hfill
\begin{subfigure}[t]{0.48\textwidth}
  \centering
  \includegraphics[width=\linewidth]{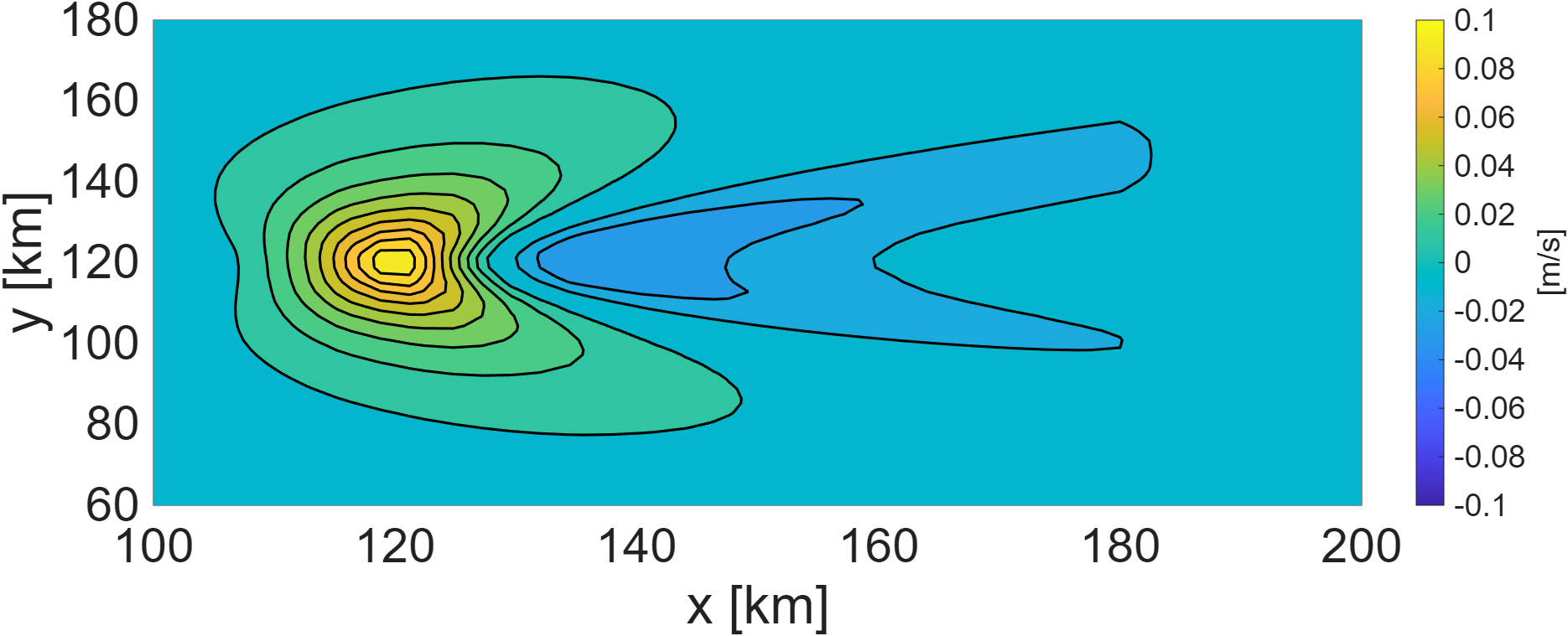}
  \caption{}
\end{subfigure}

\vspace{0.6em}

\begin{subfigure}[t]{0.48\textwidth}
  \centering
  \includegraphics[width=\linewidth]{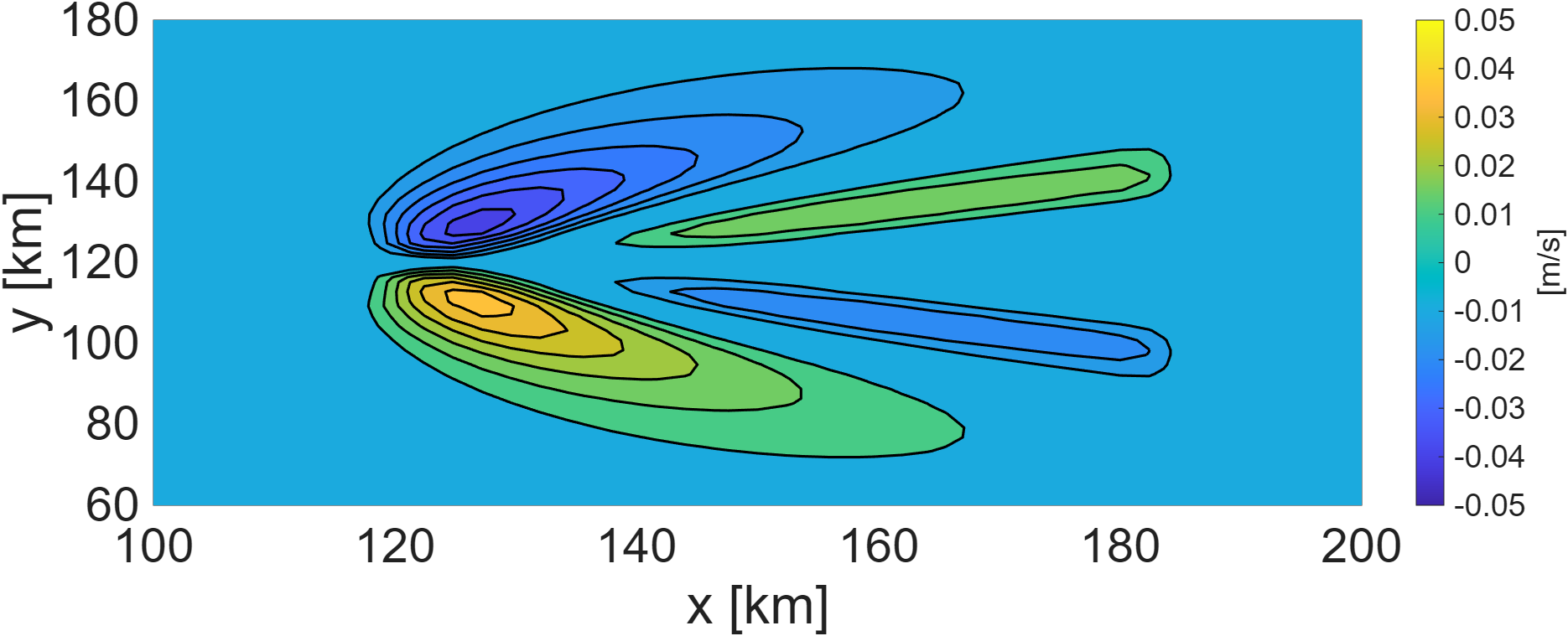}
  \caption{}
\end{subfigure}\hfill
\begin{subfigure}[t]{0.48\textwidth}
  \centering
  \includegraphics[width=\linewidth]{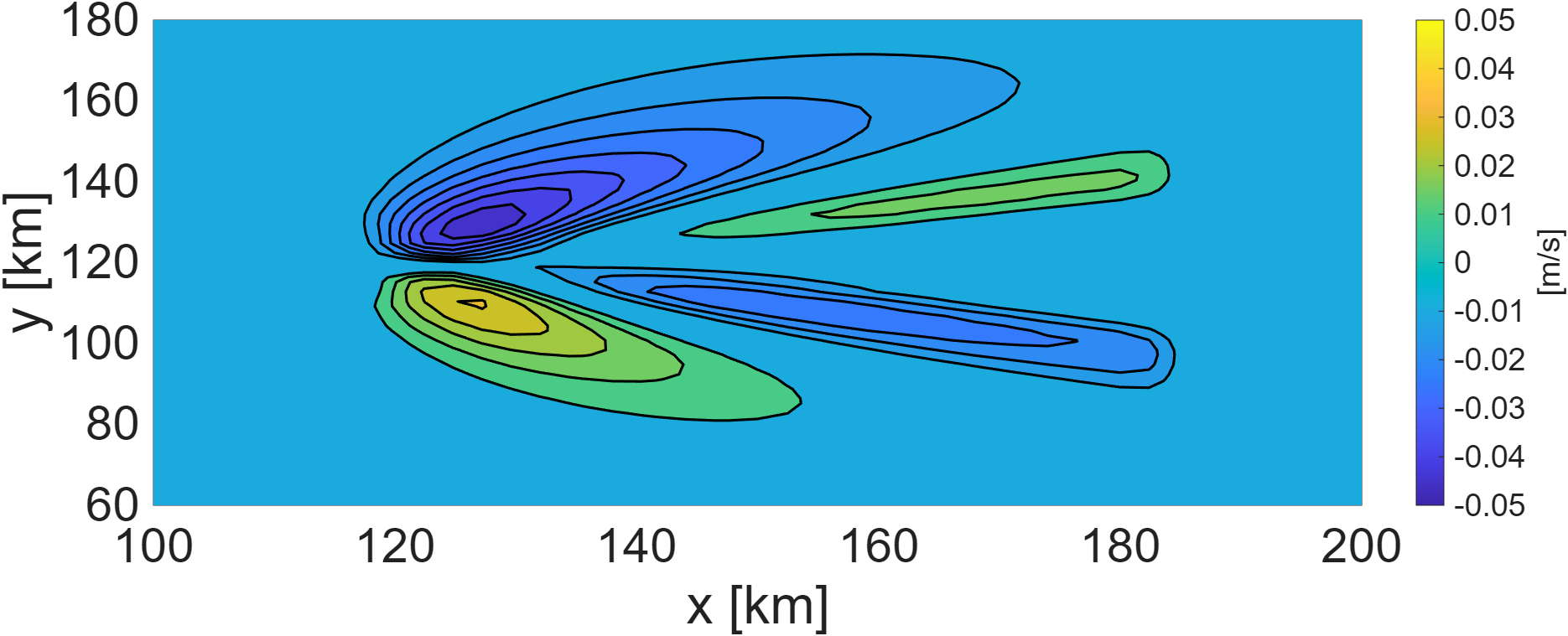}
  \caption{}
\end{subfigure}

\vspace{0.6em}

\begin{subfigure}[t]{0.48\textwidth}
  \centering
  \includegraphics[width=\linewidth]{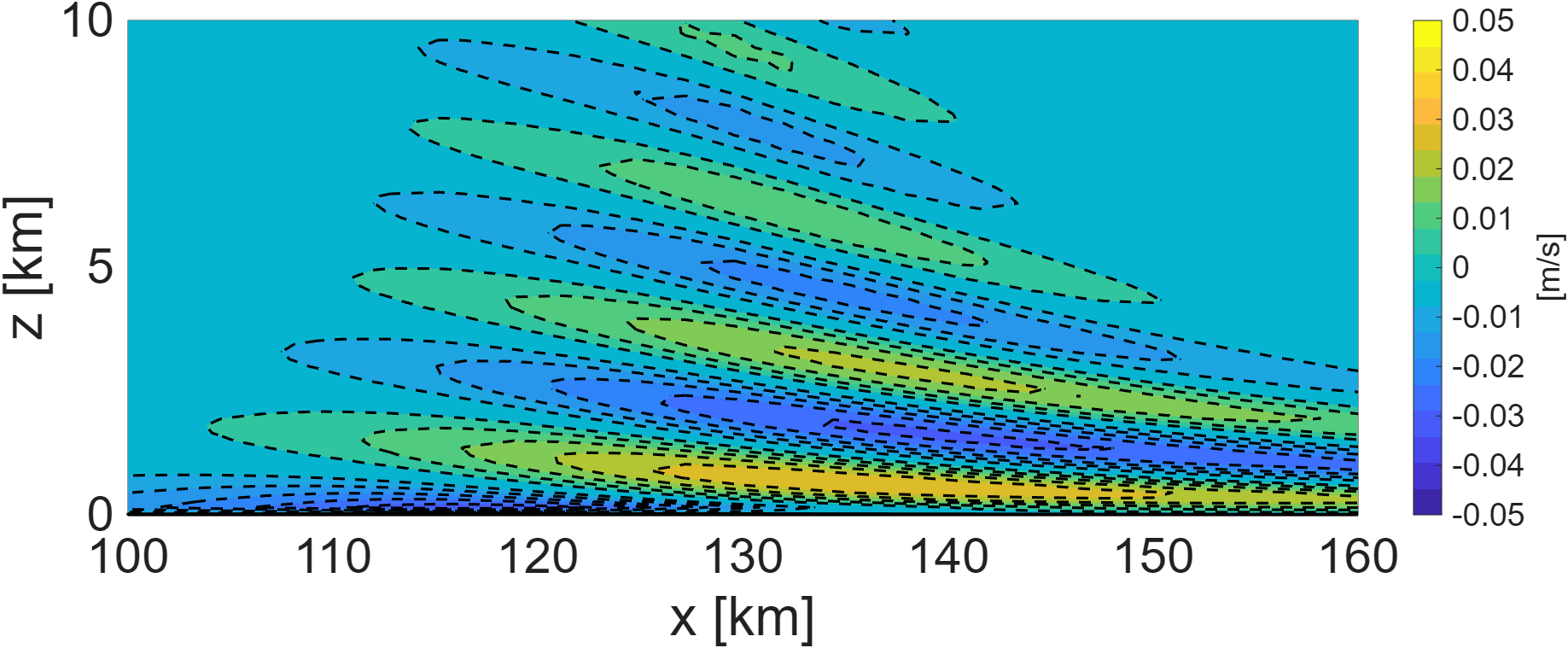}
  \caption{}
\end{subfigure}\hfill
\begin{subfigure}[t]{0.48\textwidth}
  \centering
  \includegraphics[width=\linewidth]{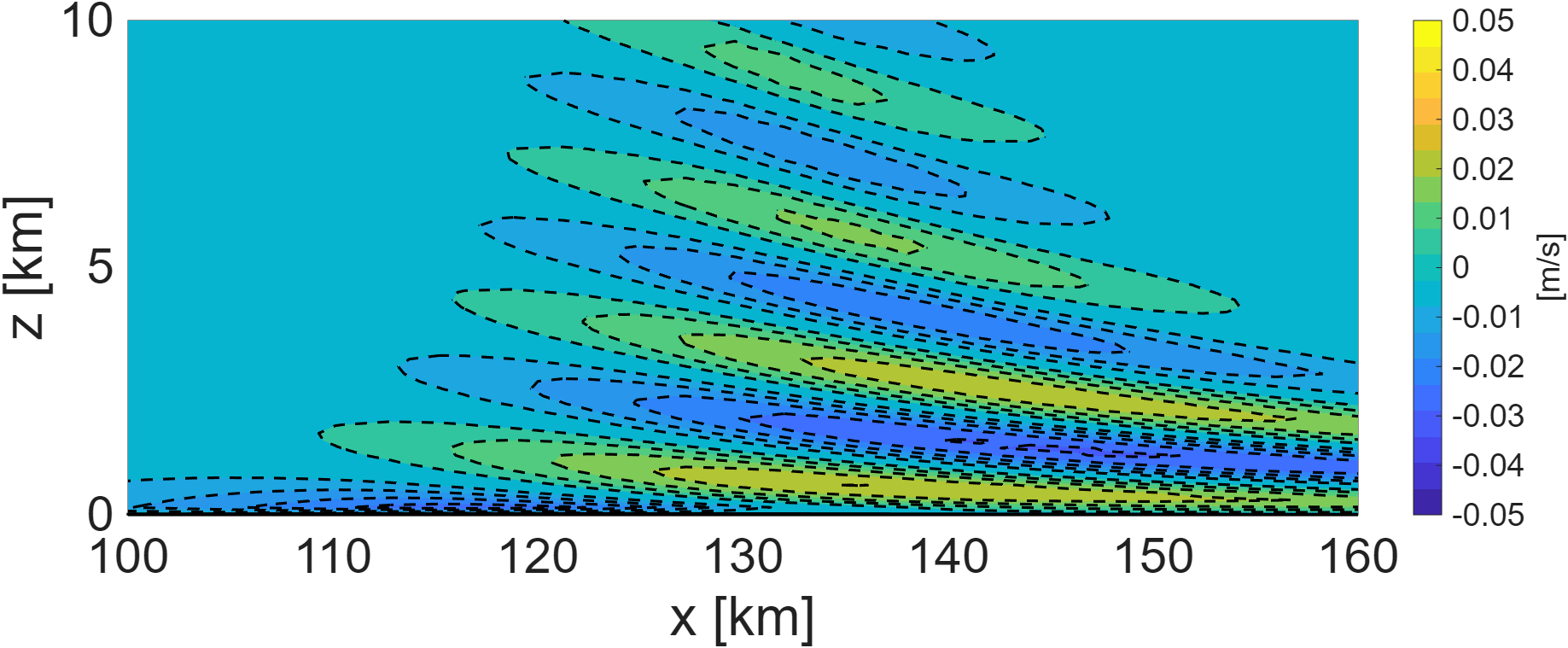}
  \caption{}
\end{subfigure}
  \begin{subfigure}[t]{0.48\textwidth}
  \centering
  \includegraphics[width=\linewidth]{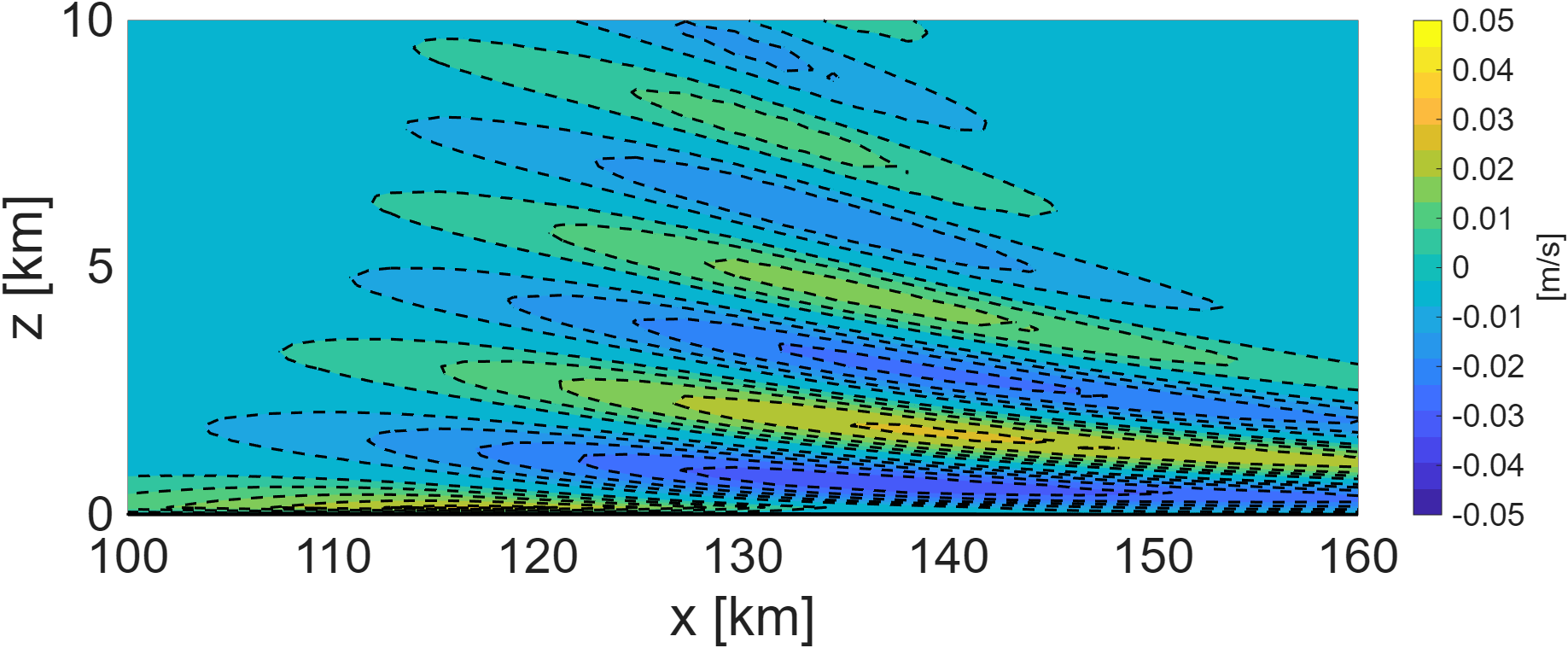}
  \caption{}
\end{subfigure}\hfill
\begin{subfigure}[t]{0.48\textwidth}
  \centering
  \includegraphics[width=\linewidth]{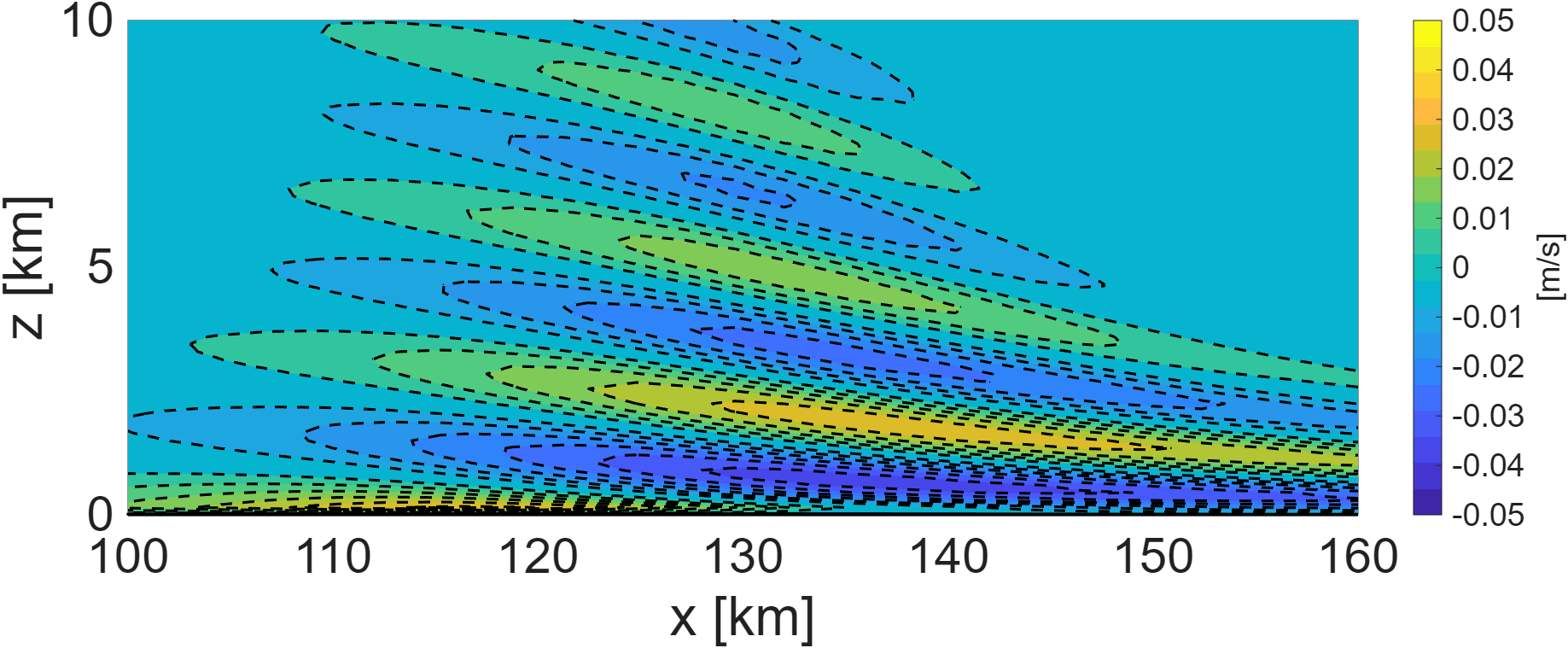}
  \caption{}
\end{subfigure}

\caption{As in Figure \ref{fig:small_domain_comparison} but for the large domain test of flow over a 3D bell-shaped hill introduced in this paper. The third and fourth row show vertical slices at $y=\SI{100}{\kilo\meter}$ and at $y=\SI{140}{\kilo\meter}$, respectively. Contours are in the range $\SI[parse-numbers=false]{[-0.1,0.1]}{\meter\per\second}$ with a $\SI{0.02}{\meter\per\second}$ interval in the top row panels, and in the range $\SI[parse-numbers=false]{[-0.05,0.05]}{\meter\per\second}$ with a $\SI{0.005}{\meter\per\second}$ interval in the remaining panels.}
\label{fig:medium_domain_comparison}
\end{figure} 

As in the smaller-domain configuration, the interaction between the stratified background flow and the orography generates a wake and a downstream gravity-wave pattern with comparable amplitudes in both the rotating and non-rotating simulations, confirming that the inclusion of the Coriolis term does not introduce spurious or unphysical effects (Figure \ref{fig:medium_domain_comparison}). However, in the present larger-scale configuration the influence of rotation is considerably more pronounced. While the non-rotating solution maintains a symmetric structure with respect to the wake axis (left panels), the rotating case develops a clear asymmetry in both the horizontal and vertical velocity fields (right panels).
This behaviour is particularly evident in the horizontal slices of the $v$ component shown in panels (c) and (d) of Figure~\ref{fig:medium_domain_comparison}, where the rotating solution exhibits a visible lateral distortion of the wave pattern. The effect becomes even clearer in the vertical $(x-z)$ sections reported in the third and fourth rows of Figure~\ref{fig:medium_domain_comparison}. These slices are taken at fixed $y$ locations and show that, at comparable downstream distances from the mountain, the non-rotating case retains nearly symmetric contour structures with opposite signs across the wake, whereas the rotating simulation displays tilted and distorted wave fronts induced by the Coriolis acceleration.
The larger horizontal extent of the domain allows the rotational effects to develop over longer spatial and temporal scales, making the associated asymmetry significantly more visible than in the smaller-domain test. At the same time, the numerical solution remains stable and free from noticeable spurious oscillations, indicating that the proposed rotating IMEX-DG formulation preserves robustness also in larger fully three-dimensional atmospheric configurations.

\begin{figure}[H]
  \centering
 \includegraphics[width=0.6\textwidth]{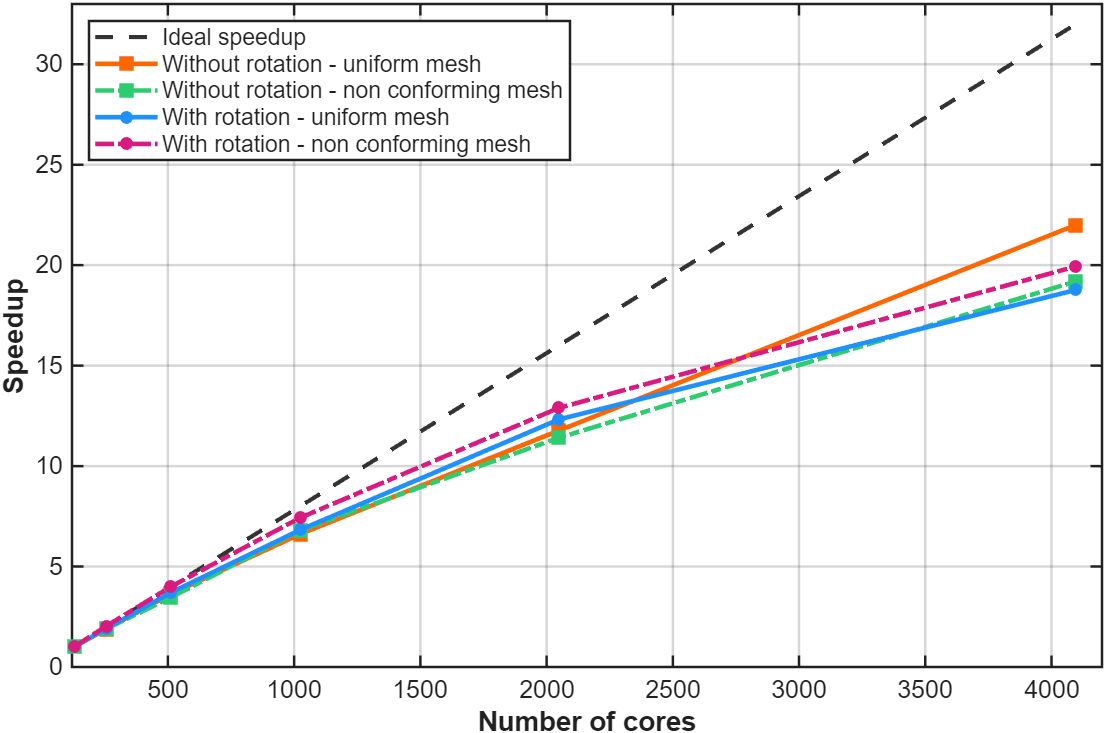}
  \caption{Large domain test: strong-scaling results for the rotating test cases using both uniform (solid lines) and non-conforming (dashed-dotted lines) meshes.} 
  \label{fig:strong_scaling_medium_domain}
\end{figure}

For this test, we also perform a strong scaling study to assess the parallel scalability of the code in a more demanding configuration, where the effects of the rotation are more pronounced. The study is carried out using both the uniform and non-conforming meshes. For the uniform case, we consider a mesh consisting of $80\times80\times60=384\,000$ elements with polynomial degree $r=4$, while for the non-conforming case, we use a mesh consisting of $183\,862$ elements. The simulations are performed using 1,2,4,8,16 and 32 full MeluXina CPU nodes, reaching 4096 cores. Figure \ref{fig:strong_scaling_medium_domain} shows the strong scaling results for the rotating and non-rotating solvers. 
For both mesh configurations, the rotating solver exhibits a scaling behaviour comparable to that of the non-rotating case over the considered range of computational resources. The speedup remains approximately linear up to 1024 cores, while the scaling efficiency gradually decreases at higher core counts. This degradation is consistent with the increasing contribution of communication and parallelisation overheads as the number of cores rises. Overall, the rotating and non-rotating configurations follow similar trends for both uniform and non-conforming meshes, indicating that the additional rotational contribution does not introduce a significant parallel-performance penalty in this more demanding configuration.

%% file: Conclusion.tex
\section{Conclusions}
\label{sec:conclusion}

We have presented an extension of the IMEX-DG framework for compressible atmospheric flows proposed in \cite{orlando_imex-dg_2023,orlando_robust_2024} to rotating configurations including the Coriolis force. The rotational contribution has been incorporated consistently within the IMEX decomposition and matrix-free DG discretisation and implemented in the \texttt{deal.II} framework. Two alternative treatments of the rotational term within the nonlinear fixed-point formulation have been introduced and compared. A relevant aspect of the proposed work is the efficient matrix-free treatment of the inverse of the operator $(\mathbf{A}+\mathbf{R})$, obtained by exploiting the discrete structure of the rotational operator.

The rotating formulation has been assessed through a sequence of two- and three-dimensional numerical experiments. The two-dimensional inertia-gravity wave benchmark provides a quantitative verification of the implementation through comparison with the analytical solution. Convergence studies show that the inclusion of the Coriolis force preserves the expected accuracy of the underlying IMEX-DG discretisation. The two alternative treatments of the rotational term lead to nearly identical error levels and convergence behaviour, indicating that both approaches provide comparable accuracy. The inertia-gravity wave configuration has been also extended to a planetary-scale domain, providing a more demanding test of the rotating formulation by considering the same dynamics over a significantly larger spatial scale.

The analysis has then been extended to fully three-dimensional flow over a bell-shaped mountain, where rotation interacts with stratification, three-dimensional dynamics, and orographic forcing. The comparison between rotating and non-rotating simulations shows that the inclusion of the Coriolis force preserves the overall wave structure and amplitude while introducing the expected asymmetry in the downstream flow field. The same configuration has been further investigated using non-conforming meshes with local refinement around the orography, demonstrating the applicability of the rotating formulation to adaptive discretisations. To the authors' knowledge, this represents the first application of the present framework to three-dimensional rotating atmospheric simulation over orography using non-conforming adaptive meshes.

Finally strong scaling tests were performed for the three-dimensional configurations to assess the computational impact of the rotational extension. The parallel performance of the rotating solver was compared with that of the non-rotating solver, showing no significant difference in strong-scaling behaviour between the two formulations. 

Overall, the results demonstrate that the proposed formulation accurately captures rotating atmospheric dynamics while retaining the accuracy, robustness, and computational efficiency of the original IMEX-DG framework. 
Future developments will focus on planetary-scale three dimensional configurations, on detailed comparisons with exactly well-balanced discretisations for longer lead times, and on extensions to spherical geometries.

\section*{CRediT authorship contribution statement}
\textbf{Letizia Bottani}: conceptualisation; investigation; formal analysis; methodology; software; validation; visualization; writing - original draft; writing - review and editing. \textbf{Tommaso Benacchio}: conceptualization; investigation; methodology; project administration; supervision; writing - review and editing. \textbf{Giuseppe Orlando}: conceptualization; investigation; methodology; software; supervision; writing - review and editing. \textbf{Luca Bonaventura}: conceptualization; investigation; funding acquisition; methodology; supervision; writing - review and editing. \textbf{Allan Engsig-Karup}: conceptualization; investigation; methodology; supervision; writing - review and editing.

\section*{Acknowledgments}

This paper expands on the first author's Master's thesis for a double degree program with Politecnico di Milano and Technical University of Denmark \cite{bottani_simulation_2026}. The first author's participation
to the 10th European Seminar on Computing 2026 was partially supported by the ESCAPE-2 project, European Union’s Horizon 2020 Research and Innovation Programme  Grant Agreement No. 800897. G.O. is part of the INdAM-GNCS National Research Group. We gratefully acknowledge Prof. Antonella Abb\`a at Politecnico di Milano for useful comments on this work.
We thank Dr. Michael Baldauf for kindly making available the code to compute the reference solution in \cite{baldauf_analytic_2013}.
We acknowledge the EuroHPC Joint Undertaking for awarding this project access to the EuroHPC supercomputer MeluXina through EuroHPC Benchmark Access calls \#EHPC-BEN-2026B01-008 and \#EHPC-BEN2026B04-091.
Generative AI tools have been used exclusively to improve language clarity, without contributing to the scientific content of the work. The authors assume responsibility for all content.

\section*{Code Availability} 

\noindent The source code used to produce the results in this study is available at:\\ \href{https://github.com/letiziabottani/Atmospheric_Flows_Rotation}
{https://github.com/letiziabottani/Atmospheric\_Flows\_Rotation}